\documentclass{amsart}

\usepackage{graphicx,array,times,mathtools}
\usepackage{color}
\usepackage[top=1.5in, bottom=1.5in, left=1in, right=1in]{geometry}
\def\squareforqed{\hbox{\rlap{$\sqcap$}$\sqcup$}}
\def\qed{\ifmmode\squareforqed\else{\unskip\nobreak\hfil
\penalty50\hskip1em\null\nobreak\hfil\squareforqed
\parfillskip=0pt\finalhyphendemerits=0\endgraf}\fi}

\newcommand{\beq}{\begin{equation}}
\newcommand{\eeq}{\end{equation}}
\newcommand{\bquo}{\begin{quote}}
\newcommand{\equo}{\end{quote}}
\newcommand{\bit}{\begin{itemize}}
\newcommand{\eit}{\end{itemize}}

\newcommand{\prt}{\partial}

\newcommand{\sub}{\subset}

\renewcommand{\phi}{\varphi}
\newcommand{\al}{\alpha}

\newcommand{\la}{\lambda}

\newcommand{\de}{\delta}

\newcommand{\be}{\beta}
\newcommand{\ga}{\gamma}

\newcommand{\Ga}{\Gamma}

\newcommand{\CQ}{{\mathcal Q}}

\newcommand{\mR}{{\mathbb{R}}}

\newcommand{\Id}{{\mathrm{Id}}}

\newcommand{\argmax}{{\mathop{\mathrm{argmax}}}}

\newcommand{\scp}[2]{{\left\langle {#1}\, , \, {#2}\right\rangle}}
\newcommand{\lform}[2]{{\left( {#1}\, |\, {#2}\right)}}

\def\Diff{\mathit{Diff}}

\newcommand{\init}{\mathit{init}}
\usepackage{pdfsync}
\graphicspath{{./Figures/}}

\title[Multiple Shape Registration]{Multiple Shape Registration using Constrained Optimal Control} 
\author[S. Arguill\`ere]{Sylvain Arguill\`ere}
\address{S. Arguill\`ere: Center for Imaging Science and Department of Applied
  Mathematics and Statistics, 
 Johns Hopkins University, 
3400 N. Charles st.
Baltimore MD 21218}
\email{sarguillere@gmail.com}
\author[E. Tr\'elat]{Emmanuel Tr\'elat}
\address{E. Tr\'elat: Sorbonne Universit\'es, UPMC Univ Paris 06, CNRS UMR 7598, Laboratoire Jacques-Louis Lions, and Institut Universitaire de France, F-75005, Paris, France.}
\email{emmanuel.trelat@upmc.fr}
\author[A. Trouv\'e]{Alain Trouv\'e}
\address{A. Trouv\'e: Ecole Normale Sup\'erieure de Cachan, Centre de
  Math\'ematiques et Leurs Applications, CMLA, 61 av. du Pdt Wilson,
  F-94235 Cachan Cedex, France}
\email{trouve@cmla.ens-cachan.fr}
\author[L. Younes]{Laurent Younes}
\address{L. Younes:
Center for Imaging Science and Department of Applied
  Mathematics and Statistics, 
 Johns Hopkins University, 
3400 N. Charles st.
Baltimore MD 21218}
\email{laurent.younes@jhu.edu}

\subjclass[2000]{58D05 49N90 49Q10 68E10}
\thanks{This work was partially supported by the ONR award N000140810606.}
\keywords{Shape analysis; optimal control; deformations; groups of diffeomorphisms.}

\begin{document}

\maketitle

\begin{abstract}
Lagrangian particle formulations  of the large deformation diffeomorphic metric mapping algorithm (LDDMM) only allow for the study of a single shape. In this paper, we introduce and discuss both a theoretical and practical setting for the simultaneous study of multiple shapes that are either stitched to one another or slide along a submanifold. The method is described within the optimal control formalism, and optimality conditions are given, together with the equations that are needed to implement augmented Lagrangian methods. Experimental results are provided for stitched and sliding surfaces.
\end{abstract}

\section{Introduction}

The large deformation diffeomorphic metric mapping (LDDMM) approach to shape matching is a powerful topology-preserving registration method with an increasing record of successful applications in medical imaging. It was first described in
\cite{joshi2000Landmark} for point sets and in
\cite{dupuis1998,trouve1998diffeomorphisms,miller2002metrics,beg2005computing}
for images and has become widely used in the medical imaging literature
and other
applications. While deeper understanding and extensions of the underlying theoretical framework was pursued \cite{miller2003metric,glaunes2004landmark,cao2005large,miller2006geodesic,younes2007jacobi,glaunes2008large,younes2010shapes,bruveris2011momentum,arguillere2014shape} and alternative numerical methods were designed \cite{avants2006lagrangian,cotter2006singular,trouve2008sparse,niethammer2011geodesic,gunther2011direct,vialard2011diffeomorphic,ashburner2011diffeomorphic,gunther2012flexible,durrleman2013sparse}, LDDMM has been applied to medical imaging data including brain \cite{miller2005increasing,zhang2006characterization,qiu2007cortical,ceritoglu2009multi,risser2011simultaneous},  heart
\cite{ardekani2009cardiac,azencott2010diffeomorphic} and  lung \cite{vidal2009template} images. This algorithm provides a non-rigid registration method between various types of objects (point sets, curves surfaces, functions, vector fields\ldots) within a unified framework driven by Grenander's concept of deformable templates \cite{grenander1993general}. It optimizes a flow of diffeomorphisms that transform an initial object (shape) into a target one.

The practical importance of shape registration is underlined by the increasing amount of work that has flourished in the literature over the past few years. LDDMM is one among
many methods that have been proposed to perform
this task. Several such methods are based on elastic matching
energies \cite{bajcsy1983computerized,droske2004variational}, and
other, like LDDMM, inspired  by viscous fluid dynamics
\cite{christensen1996deformable,thirion1998image,vercauteren2007non,ashburner2007fast,vercauteren2009diffeomorphic}.
For
surfaces, which will be our main focus, several authors have developed approaches to find
approximate conformal parametrizations with respect to the unit
disk or sphere
\cite{hurdal1999quasi,haker2000conformal,hurdal2000coordinate,gu2003genus,jin2004optimal,wang2004intrinsic,hurdal2004cortical,gu2004genus,hurdal2009discrete}.
More
recently, quasi-conformal parametrizations based on the minimization
of the Beltrami coefficient have been designed
\cite{zeng2011registration,lui2010shape,wang2009teichmuller}. Another class of
non-rigid registration methods include those based on optimal mass
transportation
\cite{haber20093d,haker2004optimal,lipman2009surface,lipman2011conformal},
while \cite{memoli2005theoretical,memoli2007use,memoli2011gromov}
introduce comparison methods based on Gromov--Hausdorff or
Gromov--Wasserstein distances. Computational methods based on integer
programming and graph optimization have also been recently introduced
\cite{zeng2010dense,zikic2010linear,glocker2008dense}. 
We also refer the reader to the survey papers \cite{zitova2003image,crum2004non,wyawahare2009image,younes2012spaces} and textbooks \cite{goshtasby2012image,younes2010shapes}  for additional entries on the literature.

In this paper, we discuss an extension of the LDDMM framework, in which multiple shapes are registered simultaneously within a deformation scheme involving contact constraints among the shapes. This is represented and solved as a constrained optimal control problem, in the spirit of the general framework recently introduced in \cite{arguillere2014shape}.

Indeed, one of the characteristics of LDDMM  is that it derives shape deformation from a global diffeomorphisms of the whole ambient space considered as a homogeneous medium, and does not allow for a differentiation of the deformation properties assumed by the shapes, or, more precisely, the objects they represent. This crude modeling may provide results that are not realistic in some applications. Consider the situation in which one studies several shapes, representing, for example,  different sub-structures of the brain. In this case, if one assumes that all shapes are deformed by a single flow of diffeomorphisms, shapes coming too close to one another will undergo a tremendous deformation, which creates artifacts that can mislead subsequent analyses. One would rather associate a different diffeomorphism to each shape,  independent from the others, but the issue is that the resulting collection of diffeomorphisms may not be consistent: the shapes could overlap along the deformation. The solution briefly introduced in \cite{arguillere2014shape} and developed in this paper is the following: embed the shapes into a "background", complement of the shapes, deformed by a new, independent deformation, and add \textit{constraints} such that, as all the shapes are simultaneously transformed, their boundary moves with the boundary of the background so that the configuration consistency is preserved. This is the approach that we develop here, focusing on surface registration. Note that a multi-diffeomorphism approach has been recently developed for image matching \cite{Risser2013182}, each diffeomorphism being restricted to a fixed region of the plane. The main (and fundamental) contrast with what we develop here is that, in our case, these subregions are variable and optimized, while they were fixed in \cite{Risser2013182}. The models along which sliding constraints are addressed in this paper and ours  also differ.  

\medskip 

This paper is organized as follows. We start by recalling the classical LDDMM algorithm in Section \ref{sec:LDDMM}, setting the definitions, notation and appropriate framework for the rest of the paper. Then, in Section \ref{sec:multi.shape}, we introduce rigorously the concept of multishape, describe identity and sliding background constraints, and describe the augmented Lagrangian algorithm for general constraints that will be used for in our numerical simulations. Section \ref{sec:discrete} follows, specializing the algorithm to the case of identity and sliding constraints in great details. Finally, Section \ref{sec:exp} applies our method to synthetic examples to to real data as well.
\\

\section{Large Deformation Diffeomorphic Metric Mapping}\label{sec:LDDMM}
\subsection{Notation}
In this paper, we define a {\em shape} as a $C^p$ embedding $q: M \to \mR^d$,
where $M$ is a compact manifold. We denote by $\mathcal M$ the
corresponding {\em shape space}, which is an open subset of the Banach space
$\mathcal Q = C^p(M, \mR^d)$.

Typical examples are as follows:
\begin{itemize}
\item $M = \{1, \ldots, m\}$ is finite, and $q$ can be identified with a
  collection $q_1, \ldots, q_m$ of distinct points in $\mR^d$.
\item $M = [0,1]$ and $q$ is a curve in $\mR^d$.
\item $M = S^{d-1}$ (the unit sphere in $\mR^d$) and $q$ is a
  hypersurface.
\end{itemize}

Our goal is to discuss models in which several shapes can deform,
while being subject to contact constraints.  The deformation process
will be similar to the one designed for {\em large deformation
diffeomorphic metric mapping} (LDDMM), which can be formulated as an {\em optimal
control problem}. Before introducing our general framework, it will be
easier to start with a description of the now well explored
single-shape problem upon which we will build. {For this, we let $(V,\Vert\cdot\Vert_V)$ be
a Hilbert space of vector fields on $\mR^d$, assumed to be continuously
imbedded in
the space $B:=C^p_0(\mR^d, \mR^d)$, the completion of the space of smooth
compactly supported vector fields for the norm $\|\cdot\|_{p, \infty}$, which
denotes the sum of supremum norms of derivatives of order $p$ or less,
with $p\geq 1$. Then $V$ possesses a {\em reproducing kernel}, that is a mapping
$K: (x,y) \mapsto K(x,y)$, defined over $\mR^d \times\mR^d$, with values
in the space of $d\times d$ matrices, such that all partial derivatives
with order less than $p$ with respect to each variable exist and 
\[
K(\cdot, y)a \in V \quad\textrm{with}\quad \scp{K(\cdot, y)a}{w}_V = a^Tw(y) ,
\] 
for all $(a,y)\in (\mR^d)^2$. The LDDMM algorithm uses flows of time-dependent vector fields $v(\cdot)\in L^2(0,1;V)$.}

\subsection{Registering Two Shapes Using LDDMM}
\label{sec:lddmm}
\subsubsection{General Problem}
The general LDDMM problem is formulated as the infinite-dimensional optimal control
problem consisting of minimizing the cost functional
\begin{equation}\label{eq:lddmm.gen}
F(v) = \frac12 \int_0^1 \|u(t)\|_V^2 \,dt + U(q(1)) ,
\end{equation}
subject to the constraint
\begin{equation}\label{eq:lddmm.gen2}
\begin{split}
\prt_t q(t)  &= u(t) \circ q(t)\quad\textrm{for a.e. $t\in [0,1]$},\\
q(0) &= q_\init.
\end{split}
\end{equation}
This differential constraint is a control system, where the control is the time-dependent vector field $u(\cdot)\in L^2(0,1;V)$, the solution of which is $q(t, \cdot) = \phi(t, q_\init(\cdot))$ where $\phi$ is the 
flow of diffeomorphisms generated by $u(\cdot)$, defined as the unique solution of the Cauchy problem $\prt_t \phi(t) = u(t) \circ \phi(t)$, $\phi(0)=\mathrm{id}_{\mR^d}$. For every time $t$, we have $\phi(t,
\cdot) \in \Diff^p$, the set of $p$-times differentiable
diffeomorphisms in $\mR^d$.

The function $U$ is a matching cost function, that is, a penalization that pushes the solution of \eqref{eq:lddmm.gen}-\eqref{eq:lddmm.gen2} towards a target. It will be assumed to be Fr\'echet differentiable from
$\CQ$ to $\mR$. To simplify the discussion, and because this covers
most of the interesting cases in practice, we will assume that there exists some fixed measure $\nu_M$ on $M$ such that its
derivative, denoted $dU(q)$ or $dU_q$ when evaluated at $q\in \mathcal Q$, can be
expressed in the form $dU_q = z_q \nu_M$ for some ($\nu_M$-measurable) $z_q: M \to \mR^d$, meaning that
\[
\forall h\in \mathcal{C}^p(M,\mR^d),\quad\lform{dU_q}{h} = \int_{M} h(m)\cdot z_q(x)\, d\nu_M(x).
\] 
Throughout the paper, for any Banach space $X$, the notation
$\lform{\mu}{v}$ will be used to designate the application $\mu(v)$ of a linear
form $\mu\in X^*$ to a vector $v\in X$.

Under these assumptions, one can prove that the gradient of the objective function $F$ defined by \eqref{eq:lddmm.gen} (which is a mapping on $L^2([0,1],V)$)
is given by 
\[
\nabla_V \tilde F(u)(t,\xi)= u(t,\xi)  - \int_M K(\xi, q(t,x))\al(t,x) \, d\nu_M(x),
\]
where $K$ is the  reproducing kernel of $V$, and $\al:[0,1]\times M\rightarrow\mR^d$ is a time-dependent function defined by $\al(1, \cdot) = - z_{q(1)}$ and
\begin{equation}
\label{eq:p.evol.lddmm}
\prt_t \al =  - (du\circ q)^T \al.
\end{equation}

This result implies, in particular, that the solutions of \eqref{eq:lddmm.gen}-\eqref{eq:lddmm.gen2} must satisfy the
Pontryagin maximum principle (see \cite{Pontryagin,Trelat}), which is the following first-order necessary condition for optimality. Let $H_u$ be the Hamiltonian defined by
\[
H_u(\rho, q) = \lform{\rho}{u\circ q} - \frac12 \|u\|^2_V,
\]
for every $u\in V$, every $q\in \mathcal Q$ and every $\rho\in \mathcal Q^*$. If $u(\cdot)$ is an optimal control, solution of the optimal control problem \eqref{eq:lddmm.gen}-\eqref{eq:lddmm.gen2}, then it must be such that
\begin{equation}\label{eq:opt.cond}
u(t) = \argmax_{w} H_w(\rho(t), q(t)) ,
\end{equation}
where $(\rho, q)$ are solutions of
\[
\begin{cases}
\prt_t q(t)  = \prt_\rho H_u(\rho(t), q(t)),\\
\prt_t \rho(t) = - \prt_q H_u(\rho(t), q(t)),
\end{cases} 
\]
($\rho$ is the so-called {\em co-state}, or {\em adjoint state})
and $\rho(1) = - dU_{q(1)}$. Indeed, it suffices to take $\rho(\cdot) =
\al(\cdot)\nu_M$ and to use properties of the reproducing kernel to check
that all the conditions are satisfied. Moreover, \eqref{eq:opt.cond} then implies that
$u=\int_M K(\cdot,q(x))\al(x)\, d\nu_M$
at every time.

\subsubsection{Examples of shapes and matching cost functions}
\label{sec:examples}

\paragraph{\bf Example 1.} To start with a simple example, let $M =\{1, \ldots, m\}$ so that a shape $q=(q(1),\dots,q(m))$ is a collection of landmarks, and consider the
landmark matching cost function defined by
$U(q) = \sum_{k=1}^m |q(k) - y_k|^2$, for fixed  $y=(y_1, \ldots, y_m)\in\mR^N$, in which $|\cdot|$ is the Euclidean norm on $\mR^d$. We then have
\[
\lform{dU_q}{h} =  2\sum_{k=1}^m (q(k) - y_k)^T h(k)\,,
\]
or, to interpret this result in the general form provided above,
$dU_q = z\nu_M$ with $z(k) = 2(q(k) - y_k)$ and $\nu_M$ the counting
measure on $\{1, \ldots, m\}$.\\

\paragraph{\bf Example 2.} If $\mu$ is a scalar measure on $\mR^d$ and $z$ a $\mu$-integrable $\mR^d$-valued
function defined on $\mathrm{support}(\mu)$, we will denote by $z\mu$
the vector measure such that
\[
\lform{z \mu}{w} = \int_{\mathrm{support}(\mu)} z\cdot w \,d\mu,
\]
where $z\cdot w $ denotes the standard euclidean dot product.

Vector measures of the form $z\mu$ are continuous linear forms over
any space that is continuously imbedded in $C^0_0(\mR^d,\mR^d)$, and in particular over
any reproducing kernel Hilbert space $W$. For such a space, equipped with a
reproducing kernel $\chi$, the operator norm of $z\mu$ is given by
\[
\|z\mu\|^2_\chi = \iint  z(x)\cdot (\chi(x,y)z(y)) \, d\mu(x) \, d\mu(y) ,
\]
and more generally, the norm of the difference between two such
measures is
\begin{multline*}
\|z\mu - \tilde{z} \tilde{\mu}\|^2_\chi = \iint  z(x)\cdot(\chi(x,y)z(y))
\, d\mu(x) \, d\mu(y) \\
-2 \iint  z(x)\cdot(\chi(x,y)\tilde{z}(y)) \, d\mu(x) \, d\tilde{\mu}(y) + \iint
\tilde{z}(x)\cdot(\chi(x,y)\tilde{z}(y)) \, d\tilde{\mu}(x) \, d\tilde{\mu}(y).
\end{multline*}
Note that $W$ and $V$ have no relationship one to each other, except that both have a continuous inclusion in $C^0_0(\mR^d,\mR^d)$, so that $\chi$ is different from $K$.

One can deduce from this the 
surface-matching cost function introduced in
\cite{vaillant2005currents}, in which an oriented surface $S$ is represented as
a geometric current and a dual-RKHS norm between currents is
used. Identifying surface currents with vector measures, this leads to the representation of $S$ given by $\mu_S = N_S \sigma_S$, where $N_S$ is the unit normal to $S$ and $\sigma_S$ its volume
form. Assume that $M$ (the parameter space) is an oriented 2D
manifold so that $S=q(M)$ is a surface, and let $\eta$ be a positively oriented volume form on
$M$. For $m\in M$ let $N_q(x)\in \mR^3$ denote the ``area-weighted
normal'' to $S = q(M)$ at $q(x)$, defined by $N_q(x) = dq_x e_1 \times
dq_x e_2$ where $(e_1, e_2)$ is an arbitrary basis of $T_xM$ such that
$\eta_x(e_1, e_2)=1$. Then 
\[
\lform{\mu_S}{w} = \int_M (w\circ q\cdot N_q) \, d\eta,
\]
for every $w\in \mathcal{C}^0_0(\mR^3,\mR^3)$.

Now, given a reproducing kernel $\chi$ and 
a target surface $\tilde{S}=\tilde{q}(M)$, we define the surface-matching cost by
\[
U(q) = \|\mu_{q(M)} - \mu_{\tilde{q}(M)}\|^2_\chi.
\]

\medskip

\paragraph{\bf Example 3.} This cost function is actually a special case of the most general framework in which one compares compact $k$-dimensional oriented submanifolds of $\mR^d$, which we briefly discuss hereafter. Given such a manifold, $S$, with a global parametrisation $q:M\to S$, one can associate to any $\omega \in C_0^p(\mR^d, (\Lambda^k\mR^d)^*)$ (the set of $C^p$ differential $k$-forms on $\mR^d$ that vanish at infinity), its integral
\[
\lform{C_S}{\omega} = \int_S \omega=\int_Mq^*\omega ,
\]
where $q^*\omega$ denotes the pull-back of $\omega$ on $M$.
An RKHS, $\tilde W$, of such forms, is a Hilbert space continuously embedded in $C^p_0(\mR^d, (\Lambda^k\mR^d)^*)$, with kernel $\tilde \chi(x,y)$ taking values in the space of bilinear functions on $\Lambda^k\mR^d \times \Lambda^k\mR^d$. The linear form $C_S$ then belongs in $\tilde W^*$, and is  a special form of a {\em geometric current}, as defined in \cite{federer1969geometric}. If $S = q(M)$ and $\tilde{S}$ is a target manifold, they can be compared using the operator norm
\begin{equation}
\label{eq:cur.met}
U(q) = \|C_{q(M)} - C_{\tilde{S}}\|_{\tilde W^*}^2.
\end{equation}


Now, if we consider $\hat{q}:\hat{M}\doteq (-1,1)\times M\to\mathbb{R}^d$ a smooth perturbation of $q$  such that $q_\varepsilon=q+\varepsilon\delta q+o(\varepsilon)$ where $q_\varepsilon(x)=\hat{q}(\varepsilon,x)$ for $\varepsilon\in (-1,1)$ and $x\in M$, we have, for $M_\varepsilon\doteq \{\varepsilon\}\times M$,
\[
\lform{dU_{q}}{\delta q} =\frac{d}{d\epsilon}\left(\int_Mq^*_\epsilon\omega\right)_{|\epsilon=0}=\frac{d}{d\epsilon}\left(\int_{M_\epsilon}\hat{q}^*\omega\right)_{|\epsilon=0}=\int_{M_0}\mathcal{L}_{\partial/\partial\varepsilon}(\hat{q}^*\omega) ,
\]
where $\mathcal{L}_{\partial/\partial\varepsilon}$ is the Lie derivative along the vector field $\frac{\partial}{\partial\varepsilon}$ on $\hat{M}$ (which is equal to $(1,0)\in \mathbb{R}\times T_mM$ at any location $(\varepsilon,x)\in\hat{M}$) and 
\begin{equation}
\label{eq:omega}
\omega = 2K_{\tilde{W}}(C_{q(M)}-C_{\tilde{S}}) ,
\end{equation}
with $K_{\tilde{W}}$ the isometry from $\tilde{W}^*$ to $\tilde{W}$.
We next show that 
\[
\lform{dU_{q}}{\delta q} = \int_M \alpha_{q}\cdot \delta q\, \mathrm{vol_{M}} + \int_{\partial M} \beta_{q}\cdot \delta q\, \mathrm{vol_{\partial M}} ,
\]
where $\mathrm{vol}_{M}$ and  $\mathrm{vol}_{\partial M}$ are the positive Riemannian volume forms on $M$ and $\partial M$, and $\alpha_{q}: M \to\mR^d$ (resp. $\beta_{q}: \partial M \to\mR^d$) is such that $\alpha_{q}(x)$ is normal to $S=q(M)$ (resp. $\beta_{q}(x)$ is normal to $\partial S=q(\partial M)$) at $q(x)$.  
Using the Cartan magic formula we get
\[
\mathcal L_{\partial/\partial\varepsilon}( q^*\omega) =  i_{\partial/\partial\varepsilon} d(q^*\omega) + d(i_{\partial/\partial\varepsilon}(q^*\omega))\,,
\]
so that, applying the Stokes theorem,
\[
\lform{dU_{q}}{\delta q} = \int_{M_0} i_{\partial/\partial\varepsilon} (\hat{q}^*d\omega) + \int_{\partial M_0} i_{\partial/\partial\varepsilon}(\hat{q}^*\omega) ,
\]
where $i_{\partial/\partial\varepsilon}$ denotes the contraction operator. Note that, for $\xi_1, \ldots, \xi_k\in T_x M$,
\begin{equation*}
  \begin{split}
    i_{\partial/\partial\varepsilon}
    (\hat{q}^*d\omega)_{(0,x)}(\xi_1,\cdots,\xi_k) & =
    d\omega_{q(x)}(\delta q(x), dq_x\xi_1, \ldots,
    dq_x\xi_k) \\
&= d\omega_{q(x)}(\delta q^\perp(x), dq_x\xi_1,
    \ldots, dq_x\xi_k) ,
  \end{split}
\end{equation*}
where $\delta q^\perp(x)$ denotes the projection of $\delta q(x)$ on $T_{q(x)} S^\perp$  (since the form vanishes if $\delta q(x)\in T_{q(x)} S=T_{q(x)}q(M)$). Since the set of $k$-forms is a one-dimensional space on $M$, this means that we can write, for some function $\alpha_{q}$ such that  $\alpha_{q}(x) \in (T_{q(x)} S)^\perp$, with $x\in M$,
\[
i_{\partial/\partial\varepsilon}(q^*d\omega)_{(0,x)} = ((\delta q\cdot \alpha_{q})\,  \mathrm{vol}_{M})_x\,.
\]
Similarly,
\[
i_{\partial/\partial\varepsilon}(q^*\omega_{q})_{(0,x)} = ((\delta q\cdot \beta_{q}) \mathrm{vol}_{\partial M})_x ,
\]
for $\beta_{q}(x) \in (T_{q(x)}\partial S)^\perp$. 

The data term defined in \eqref{eq:cur.met} is derived from this general construction, using the fact that two-forms in $\mR^3$ (or $(d-1)$-forms in $\mR^d$) can be identified with vector fields via $\omega_x(e_1,e_2) = w(x)\cdot(e_1\times e_2)$. The current $C_S$ is then identified with the vector measure $\mu_S$. The form $\omega$ introduced in \eqref{eq:omega} becomes, introducing a global parametrizaion $\tilde{q}:M\to \tilde{S}$ of the target $\tilde{S}$, the vector field
\[
w(\cdot) = 2\int_{\mR^d} \chi(\cdot, y) \, d(\mu_{q(M)}-\mu_{\tilde{q}(M)})(y).
\]
With this identification, we have $\alpha_{q} = \mathrm{div}(w)\circ q\, N_q$ where $N_q$ is the oriented ``area-weighted'' normal to $S$ at $q(x)$ defined previously, and $\beta_{q} = \tau_q \times w\circ q$, where $\tau_q$ is the oriented ``length-weighted'' tangent to $q(\partial M)$ given as $\tau_q(x)=dq_xe_1$, where $e_1$ is the unit positively oriented tangent vector at $x$ along $\partial M$.

\medskip

\paragraph{\bf Example 4.} Returning to surfaces,  the discrete case, in which triangulated
surfaces are compared, is, for practical purposes, even more important. We here also consider
the case $M = \{1, \ldots, m\}$, with an additional family $F$ of
facets, which are
ordered triples $(i,j,k)$ with $i,j,k \in M$. (We assume that $F$ is a
consistent with a manifold structure: The set $V_i$ of indices that
share a facet with $i$ must form a chain and no pair of indices can be
included in more than two facets.) 

Given a one-to-one mapping $q: M \to \mR^3$,
define $S_q$ as the collection
of triangles $S_q = \{(q(i), q(j), q(k)), (i,j,k)\in F\}$. If $f =
(i,j,k)$, let $q(f) =  (q(i), q(j), q(k))$, $N_q(f) = (q(j) - q(i))
\times (q(k)-q(i))$ and $c_q(f) = (q(i)+q(j)+q(k))/3$ respectively denote
the triangle, area-weighted normal and center associated to the facet
$f$. Following \cite{vaillant2005currents}, we define the vector
measure associated to $q$ by
\[
\mu_q = \sum_{f\in F} N_q(f) \de_{c_q(f)}.
\]
Here, $\de_x$ (with $x\in\mR^3$) denotes the atomic measure of mass 1 with support $\{x\}.$ The (discrete) surface matching cost associated to a target $\tilde q$
is then defined by
\[
U(q) = \|\mu_q - \mu_{\tilde q}\|_\chi^2.
\] 
Note that $\tilde q$ does not need to be consistent with $q$, and can be
defined on a different set of indices, $\tilde M= \{1, \ldots, \tilde
m\}$ and triangle structure $\tilde F$.
One then has
$dU_q = \al_{q} \nu_M$, where, as above, $\nu_M$ is the counting measure on $M$, and
\[
\al_{q}(i) = \sum_{f\in F, i\in f} (dZ_{c_q(f)}^TN_q(f)/3 + e_{q}(f,i) \times Z(c_q(f))) ,
\]
with $e_{q}(f,i)= q(k)-q(j)$ the oriented edge opposed to $q(i)$ in $q(f)$, and
\[
Z(\cdot) = \sum_{f\in F} \chi(\cdot, c_q(f)) N_q(f) - \sum_{\tilde f\in
  \tilde F} \chi(\cdot, c_{\tilde q}(\tilde f) N_{\tilde q}(\tilde f).
\]

\subsubsection{Reduced Problem}
Since the optimal control must satisfy 
\begin{equation}
\label{eq:v.al}
v = \int_M K(\cdot, q(t,x))\al(t,x) \,d\nu_M
\end{equation}
for some function $\al$ defined on $M$, it is natural to 
parametrize $v$ by $\al$ and use this function as a new control. 
We define the inner product
\[
\scp{\al}{\be}_q = \int_{M\times M} \al(x)^T K(q(x), q(\tilde x)) \be(\tilde
x)\, d\nu_M(x)\, d\nu_M(\tilde x) 
\]
between two measurable functions $\al$ and $\be$ defined on $M$.
If $v$ is given by \eqref{eq:v.al}, 
the reproducing property of the kernel implies that $\|v\|_V^2 = \|\al\|_q^2$. 
The optimal
control problem \eqref{eq:lddmm.gen}-\eqref{eq:lddmm.gen2} is then equivalent to the reduced problem consisting of minimizing the cost functional
\begin{equation}
\label{eq:lddmm.red}
F(\al) = \frac12 \int_0^1 \|\al(t)\|_{q(t)}^2 \, dt + U(q(1)) ,
\end{equation}
subject to the constraint (control system)
\begin{equation}
\label{eq:lddmm.red.q}
\prt_t q(t,x) = \int_{M} K(q(t,x), q(t,\tilde x)) \al(t,\tilde
x) d\nu_M(\tilde x),
\end{equation}
almost everywhere over the time interval $[0,1]$.

According to \cite{azencott2010diffeomorphic, Cotter:2009, vaillant2005currents,younes2009evolutions}, we have $\nabla F(\al) = \al - p$,
where $p$ (the co-state) is a time-dependent vector-valued measurable function on $M$ such that
$p(1)\nu_M = - dU_{q(1)}$ and 
\[
\prt_t p(t) = -\prt_q \left( \scp{p(t)}{\al(t)}_q -\|\al(t)\|_{q(t)}^2/2 \right) ,
\]
where $q$ is defined by \eqref{eq:lddmm.red.q}. Here, the gradient is
computed with respect to the inner product $\scp{\cdot}{\cdot}_q$.

\section{Multiple Shape problems}
\label{sec:multi.shape}
\subsection{Motivating Examples}
In the previous formulation,  the shape evolution was controlled by a
single, smooth vector field $v$, 
inducing a single diffeomorphism of $\mR^d$ restricted
to the considered shape. This approach has been successfully used to
model variations of single, homogeneous shapes, and led to important
applications in computational anatomy, including, among many other
examples, the impact of pathologies like Huntington disease
\cite{younes2012regionally}, schizophrenia
\cite{qiu2009neuroanatomical}, and Alzheimer's disease
\cite{wang2007large,durrleman2009spatiotemporal,risser2011simultaneous} on brain structures. 
This deformation model, however, is not well adapted in situations in
which several shapes interact, or situations in which shapes have
heterogeneous parts. Let us review some motivating examples.

\begin{enumerate}
\item[(1)] Consider a schematic representation of a kite, or a manta ray,
 composed  with a two-dimensional surface, representing the body, and an
  open curve attached to it representing the tail. When comparing two such
  objects, the body is assumed to only show small differences in
  shape, while the tail can vary widely.
\item[(2)] Consider a two-dimensional representation of a mouth, with two
  curves representing the upper and lower lip. Because the mouth can
  be wide open or closed, it is not possible to consider its
  deformations as resulting from the restriction of a smooth
  diffeomorphism of $\mR^2$.
\item[(3)] Finally, it is natural, when analyzing multiple organs in the
  human body, to consider multiple shapes, each of them being
  relatively stable (only subject to small deformations) while their
  position with respect to each other is subject to larger variations,
  so that the background (the intersection of their complements) is
  subject to very large deformations. Here again, modeling the whole
  process with a single diffeomorphism is not adequate. 
\end{enumerate}
\medskip

These examples suggest using multiple deformations 
applied to each  component of the considered model. Generalizing \eqref{eq:lddmm.gen}-\eqref{eq:lddmm.gen2},
consider parameter spaces $M_1, \ldots, M_n$ for an $n$-component
model. Each shape, or component, is a mapping $q^{(k)}\in \mathcal Q_k: M_k \to
\mR^d$. The shape space will then be $\mathcal Q = \mathcal Q_1
\times\cdots \times
\mathcal Q_n$. To each shape, associate a control $u_k \in V_k$, where
$V_k$ is an RKHS embedded in $C^p_0(\mR^d, \mR^d)$ with the state
evolution equation $\partial_t q^{(k)} = u_k \circ q^{(k)}$. We can then choose each $V_k$ according to how ``wildly" we want to allow the $k$-th shape to deform. These evolutions, however,
must be consistent with each other, implying contact constraints that we will consider in two forms:

\begin{itemize}
\item {\bf Identity constraints:} These are constraints that make a subset of the $k$-th shape stay stitched to a subset of the $l$-th shape, so that these subsets coincide in $\mR^d$ and move identically along the deformation. 
Given some pair $(k,l) \in \{1,\ldots, n\}^2$, and given a one-to-one mapping $g_{kl}: A_{kl}\sub M_k \to A_{lk} = g_{kl}(A_{kl})\sub M_l$, one has $q^{(k)}(x) = q^{(l)}(g_{kl}(x))$ for every $x\in A_{kl}$. 

\item {\bf Sliding constraints:} {These are constraints that force a closed submanifold of the $k$-th shape to slide on a corresponding submanifold of the $l$-th shape along the deformation, for all $(k,l)$. Here, we assume that all parameter spaces are orientable differential manifolds, and that all $q^{(k)}$'s are immersions. Given some pair $(k,l) \in \{1,
  \ldots, n\}$, and a closed submanifold without boundary $A_{kl}\sub M_k $, there exists a diffeomorphism $g_{kl}: A_{kl} \to A_{lk}\sub M_l$ onto a fixed closed submanifold $A_{lk}$ of $M_l$ such that $q^{(k)}(x) = q^{(l)}(g_{kl}(x))$ for every $x\in A_{kl}$.}
\end{itemize}

Let us turn back to the examples mentioned at the beginning of the section.
For Example (1), we can take $M_1 = S^2$ and $M_2 = [0,1]$, and, letting $x_0$ represent the north pole in $S^2$, impose $q^{(1)}(x_0) =
q^{(2)}(0)$. We can then assign different deformation models to $q^{(1)}$ and $q^{(2)}$ via
the metrics on $V_1$ and $V_2$. 

Example (2) requires a slightly more complex
construction, that only imperfectly addresses the issue. Let $M_1 = M_2 = [0,1]$ and $M_3 = \{1, 2\} \times
[0,1]$. Let $q^{(1)}$ represent the upper lip, $q^{(2)}$ the lower one and $q^{(3)}$ their union. We use the identity constraints $q^{(1)}(1) = q^{(2)}(0)$ and $q^{(2)}(1) =
q^{(1)}(0)$ for the extremities of each lip, and $q^{(3)}(1, \cdot) =
q^{(1)}(\cdot)$, $q^{(3)}(2, \cdot)= q^{(2)}(\cdot)$. We take $V_1=V_2$ and choose $V_3$
such that the latter allows for large deformations at a small cost.
With this model, it becomes easier to ``almost'' close the mouth,
although the deformation inside the mouth must remain diffeomorphic, so that the closing cannot go all the way. 

For Example (3), there are $n$ shapes, $n-1$ of which are associated with the organs, and the last of which represents the background. For example, we can take $M_k = S^2$ for $k=1,
\ldots, n-1$, and $M_n = \{1, \ldots, n-1\} \times S^2$. Assuming that
the shapes do not intersect, we can define identity or sliding constraints for
the background, enforcing $q^{(n)}(k,\cdot)=q^{(k)}(\cdot)$ or $q^{(n)}(\{k\}\times S^2) = q^{(k)}(S^2)=S_k$ for $k\in \{1, \ldots, n-1\}$ during the deformation.

\subsection{Induced Constraints}
The previous constraints can be reformulated as equality constraints involving the state and control.
Identity constraints $q^{(k)}(x) = q^{(l)}(g(x))$ are equivalent (taking time derivatives) to $u_k(t, q^{(k)}(t, x)) = u_l(t, q^{(l)}(t, g(x)))$ as soon as the constraints are satisfied at time $t=0$, which we  obviously assume. 

Making the same assumption, sliding constraints can be expressed as 
\begin{equation}
\label{eq:sliding}
 N^{(k)}(t, q^{(k)}(t, x))^T (u_k(t, q^{(k)}(t, x)) - u_l(t, q^{(k)}(t, x))) = 0,
\end{equation}
{where $N^{(k)}(t, q^{(k)})$ is a $d\times (d-\mathrm{dim}(A_{kl}))$ matrix consisting of independent vectors perpendicular to $T_{q^{(k)}}B_{kl}(t)$ (e.g, a normal frame to $B_{kl}$), with $B_{kl}=q^{(k)}(A_{kl})$ for every $(k,l)$. Let us briefly justify this statement.}

We express the sliding constraint as
$q^{(k)}(t, x) = q^{(l)}(t, g(t,x))$ for some {diffeomorphism $g(t, \cdot): A_{kl}\to A_{lk}$,} assuming a differentiable dependency on time. Taking time derivatives, we get
$$
u_k(t,q^{(k)}(t,x))=u_l(t,q^{(l)}(t,g(t,x)))+dq^{(l)}(t, g(t,x)) \partial_t g(t, x),\quad x\in M_k.
$$
Since $q^{(l)}(t,g(t,x))=q^{(k)}(t,x)$, we obtain
\begin{eqnarray*}
u_k(t, q^{(k)}(t, x)) - u_l(t, q^{(k)}(t, x))& =& dq^{(l)}(t, g(t,x)) \partial_t g(t, x)\\
&=& dq^{(k)}(t,x) dg(t,x)^{-1} \partial_t g(t,x) ,
\end{eqnarray*}
which is tangent to {$B_{kl}$ at $x$. Note that, since the image of $g(t,\cdot)$ is $A_{lk}$ for every time $t$, we do have $\partial_t g(t, x)\in T_{g(t,x)}A_{lk}=dg(t,x)(T_xA_{kl}),$ so $dg(t,x)^{-1} \partial_t g(t,x) $ is well-defined.}

{Conversely, assume that \eqref{eq:sliding} holds for every $x\in A_{kl}$, with $q^{(k)}(0, x) = q^{(l)}(0, g_0(x))$ for some diffeomorphism $g_0: A_{kl} \to A_{lk} \sub M_l$.  
Then for every time $t$, the mapping
\[
w:x\in A_{kl}\mapsto dq^{(k)}(t, x)^{-1} (\underset{\in\ T_{q^{(k)}(x)}B_{kl}}{\underbrace{u_k(t, q^{(k)}(t,x)) - u_l(t, q^{(k)}(t,x))}})\ \in T_xA_{kl}
\]
defines a time-dependent vector field on $A_{kl}$. Since $A_{kl}$ is a closed manifold, this vector field is complete, and we denote its flow by $h(t, \cdot):A_{kl}\rightarrow A_{kl}$. Then
\begin{eqnarray*}
\partial_t q^{(k)}(t,h(t,x)) &=& u_k(t, q^{(k)}(t,h(t,x))) +   dq^{(k)}(t,h(t,x)) \partial_t h(t,x) \\
&=& u_l(q^{(k)}(t,h(t,x))) ,
\end{eqnarray*}
where the last identity holds for every $x\in A_{kl}$,
so that $q^{(k)}(t,h(t, x))$ and $q^{(l)}(t, g_0(x))$ satisfy the same differential equation with the same initial condition and therefore coincide. Hence, letting $g(t,x)=g_0(h^{-1}(t,x))$, we obtain $q^{(l)}(t, g(t,x))=q^{(k)}(t,x)$ for every $x\in A_{kl}$.}

{It is possible to extend this construction to the case where $A_{kl}$ is a compact manifold with boundary. In this case, the matrix $N^{(k)}(t)$ must consist of a normal frame along $\partial B_{kl}(t)$ and possesses therefore an extra column, and we get an additional constraint along the boundary of $A_{kl}$.}

\medskip

We will need the constraints to depend smoothly on $q$, and therefore we will need a smooth representation of the normal space to $q$ (a smooth map $q\mapsto N(q)$) in order to be able to use \eqref{eq:sliding}. When this is not possible (or convenient), one can also use the alternative approach of introducing a new state, say $N^{(k)}(t, x)$, evolving according to
\begin{equation}
\label{eq:normal.evol}
\partial_t N^{(k)} = - du_k(q^{(k)})^T N^{(k)} ,
\end{equation}
which ensures that $N^{(k)}(t,x)$ remains perpendicular to $T_{q^{(k)}(t,x)}S_k$ as soon as this holds  true at $t=0$. The constraint $ N^{(k)}(t, x)^T (u_k(t, q^{(k)}(t, x)) - u_l(t, q^{(k)}(t, x)))$ is now a smooth function of the extended state.

\medskip
 
The two problems that we consider are therefore special cases of the general problem considered in \cite{arguillere2014shape}, which is the problem of minimizing the cost functional
\begin{equation}\label{eq:mult.oc}
\frac{1}{2}\sum_{k=1}^n \int_0^1 \|u_k(t)\|_{V_k}^2 \, dt + \sum_{k=1}^n U_k(q^{(k)}(1)) ,
\end{equation}
subject to the constraints
\begin{equation}\label{eq:mult.oc1}
\partial_t q^{(k)}(t) = u_k (t,q^{(k)}(t)),\quad\textrm{and}\quad C(q(t))u(t) = 0,
\end{equation}
almost everywhere over the time interval $[0,1]$, where $C: \mathcal M \to L(V, \mathcal Y)$ takes values in the space of bounded linear operators from $V$ to a Banach space $\mathcal Y$. Here, we have $V = V_1\times \cdots \times V_n$ and $q = (q^{(1)}, \ldots, q^{(n)})$.

The study of this constrained optimal control problem, and in particular, the derivation of its first-order optimality conditions (of the type of Pontryagin maximum principle), is challenging in this infinite-dimensional setting. In \cite{arguillere2014shape}, it is proved that, under some differentiability conditions, and {\em under the important assumption that $C(q)$ is surjective for every $q\in \mathcal M$}, optimal solutions must be such that there exist $p=(p^{(1)}, \ldots, p^{(n)}) \in H^1([0,1], \mathcal Q^*)$ and $\lambda\in L^2([0,1], \mathcal Y^*)$ that satisfy
\begin{equation}
\label{eq:pmp.gen}
\begin{dcases}
\partial_t q^{(k)} = u_k(q^{(k)}) , \\
\partial_t p^{(k)} = -\partial_{q^{(k)}} \lform{p^{(k)}}{u_k(q^{(k)})} - \partial_{q^{(k)}} \lform{\lambda}{C(q)u} , \\
\langle u_k,v\rangle_{V_k} = - \lform{p^{(k)}}{v\circ q^{(k)}} - \lform{\lambda}{C_k(q)v},\quad v\in V_k , \\
\sum_{k=1}^{n} C_{k}(q) u_{k} =0 ,
\end{dcases}
\end{equation}
where $C_k(q)u_k = C(q)(0, \ldots,0, u_k, 0, \ldots, 0)$.

Unfortunately, the constraints $C(q)$ that correspond to our identity or contact constraints are, in general, not surjective, and the results of \cite{arguillere2014shape} cannot be applied in a fully general infinite-dimensional context. However, surjectivity becomes almost straightforward when these constraints are discretized to a finite number. They are true as soon as the points involved in the constraints are all distinct, which is a mild assumption. We now proceed to the description of a discrete version of this approach.

\section{Discrete Approximations}\label{sec:discrete}

\subsection{Augmented Lagrangian}
As an example, and to simplify the presentation, we detail our implementation for multi-shape problems in which shapes interact (through constraints) with a background, but not directly with each other. Direct interactions between shapes can be handled in a similar way.  Our constrained optimization method uses the augmented Lagrangian method (see, e.g., \cite{nocedal2006numerical}). In a nutshell, in order to minimize a function $u \mapsto F(u)$ subject to multi-dimensional equality constraints $C(u) = 0$, the augmented Lagrangian method consists of considering the functional 
\[
L(u) = F(u) - \lambda\cdot C(u) + \frac \mu 2 |C(u)|^{2} ,
\]
in which $\lambda$ lives in the dual space of the space of constraints $Y$, and $\mu$ is a positive real number. Each iteration of the algorithm consists in minimizing $L$ with fixed $\lambda$ and $\mu$ (our implementation using nonlinear conjugate gradient) until the gradient norm passes below some upper bound, and then in updating $\lambda$ according to the rule
\[
\lambda \leftarrow \lambda -\mu C(u) ,
\]
before running a new minimization of $L$. 
The constant $\mu$ is increased only if needed, i.e., if the norm of the constraint did not decrease enough during the minimization. More details can be found in \cite{nocedal2006numerical}.

We first apply this to identity constraints, which only require the shapes to be discretized into a sets of points. We will then discuss sliding constraints, which will require more structure in order to define normal frames to the boundary. 

\subsection{Identity Constraints}

We consider $n-1$ objects, discretized into point sets, so that $M_{k}$ is a finite set of indices for each $k$. Let  $x^{(k)}_j = q^{(k)}(j)$ and  $x^{(k)} = (x^{(k)}_1, \ldots, x^{(k)}_{m_k})$, for $k=1, \ldots, n-1$, with $m_k = |M_k|$. We add as $n$-th object the background, defined on $M_n = (\{1\}\times M_1) \cup \cdots \cup (\{n-1\}\times M_{n-1})$. We let $z^{(k)}_j = q^{(n)}(k,j)$,   $z^{(k)} = (z^{(k)}_1, \ldots, z^{(k)}_{m_k})$ and
$z = (z^{(1)}, \ldots, z^{(n-1)})$  (a collection of $m = m_n = m_1 + \ldots + m_{n-1}$ points). 

Assume that end-point cost functions $U_1(x^{(1)}), \ldots, U_{n-1}(x^{(n-1)})$ are defined, typically measuring the discrepancy between each collection of points and an associated target. We assume similar functions $\tilde U_1(z^{(1)}), \ldots, \tilde U_{n-1}(z^{(n-1)})$ for the background, typically using $U_{j} = \tilde U_{j}$. The associated constrained optimal control problem consists in minimizing the cost functional
\[
\frac12 \sum_{k=1}^{n} \int_0^1 \|u_{k}(t)\|^{2}_{V_{k}} \, dt 
+ \sum_{k=1}^{n-1} U_k(x^{(k)}(1)) +  \sum_{k=1}^{n-1} \tilde U_k(z^{(k)}(1)) ,
\]
subject to the constraints (almost everywhere along $[0,1]$)
\[
\begin{cases}
\prt_t x^{(k)}_{j} = u_{k}(x^{(k)}_{j})\quad j=1, \ldots, m_{k},\  k=1, \ldots, n-1 ,\\
\prt_t z^{(k)}_{j} = u_{n}(z^{(k)}_{j}), \quad j=1, \ldots, m_{k},\  k=1, \ldots, n-1  ,\\
z^{(k)} = x^{(k)}, \quad k=1, \ldots, n-1.
\end{cases}
\]

For $y$ and $y'$ ordered families of points in $\mathbb R^{d}$, let $K^{(k)}(y, y')$ be the matrix formed with all $d\times d$ blocks $K_{V_k}(y_i, y'_j)$, and let $K^{(k)}(y) = K^{(k)}(y,y)$, for $k=1, \ldots, n$, where $K_{V_k}$ is the kernel of $V_k$. Since the problems only depend on the values taken by $u_{1}, \ldots, u_{n}$ on their corresponding point set trajectories $x^{(1)}, \ldots, x^{(n-1)}, z$, the optimal vector fields take the form
\begin{equation*}
\begin{split}
u_k(\cdot) &= K^{(k)}(\cdot, x^{(k)})\alpha^{(k)},\quad k=1, \ldots, n-1, \\
u_n(\cdot) &= K^{(n)}(\cdot, z)\beta,
\end{split}
\end{equation*}
for some families of $d$-dimensional vectors $\alpha^{(1)}, \ldots, \alpha^{(n-1)}, \beta$. The problem can therefore be reduced to the finite-dimensional optimal control problem consisting in minimizing the cost functional
\begin{multline*}
E(\al, \be, x, z) = \frac12 \sum_{k=1}^{n-1} \int_0^1 \al^{(k)} \cdot (K^{(k)}(x^{(k)})\al^{(k)}) \, dt + \frac12 \int_0^1 \be\cdot (K^{(n)}(z)\be) \, dt \\
+ \sum_{k=1}^{n-1} U_k(x^{(k)}(1)) + \sum_{k=1}^{n-1} \tilde U_k(z^{(k)}(1))
\end{multline*}
subject to the constraints (almost everywhere along $[0,1]$)  
\[
\begin{cases}
\prt_t x^{(k)} = K^{(k)}(x^{(k)})\al^{(k)} ,\\ 
\prt_t z = K^{(n)}(z) \be ,\\
z^{(k)} = x^{(k)}, \quad k=1, \ldots, n-1 .
\end{cases}
\]

Extending $E$ with the augmented Lagrangian method, we introduce coefficients $\la^{(k)}$, $k=1, \ldots, n-1$ (where $\la^{(k)}$ has the same dimension as $x^{(k)}$) and $\mu>0$, defining
\begin{multline*}
L(\al, \be, x, z) = \frac12\sum_{k=1}^{n-1} \int_0^1 \al^{(k)} \cdot (K^{(k)}(x^{(k)})\al^{(k)}) \,dt + \frac12 \int_0^1 \be\cdot (K^{(n)}(z)\be) \,dt +  \sum_{k=1}^{n-1} U_k(x^{(k)}(1))\\
+  \sum_{k=1}^{n-1} \tilde U_k(z^{(k)}(1)) - \sum_{k=1}^{n-1} \int_0^1 \la^{(k)}\cdot(x^{(k)} - z^{(k)}) \, dt  + \frac{\mu}{2}  \sum_{k=1}^{n-1} \int_0^1 |x^{(k)} - z^{(k)}|^2 \,dt ,
\end{multline*}
which will be minimized subject to the constraints (almost everywhere along $[0,1]$) 
\[
\begin{cases}
\prt_t x^{(k)} = K^{(k)}(x^{(k)})\al^{(k)} , \\
\prt_t z = K^{(n)}(z) \be .
\end{cases}
\]

From the constraints, $L$ can be considered as a function of $\alpha$ and $\beta$ only, and its differential with respect to these variables can be computed via the adjoint method as follows. Denoting the co-states by $p^{x,k}$, $k=1, \ldots, n-1$, and $p^z$, the associated Hamiltonian is 
\[
H = \sum_{k=1}^{n-1} \int_0^1 p^{x,k} \cdot K^{(k)}(x^{(k)})(\al^{(k)}) \, dt + \int_0^1 p^z \cdot K^{(n)}(z)\be \,dt - L.
\] \\
The computation of the gradient of $L$ follows the same general scheme as the one described   in Section \ref{sec:lddmm} for the basic LDDMM algorithm. Given $\alpha$ and $\beta$ and the associated trajectories $x$ and $z$, one has solve the adjoint equations
\begin{equation*}
\begin{split}
\partial_t p^{x,k} &= - \partial_{x^{(k)}} H,\quad p^{x,k}(1) = -\nabla U_k(x^{(k)}(1)), \quad k=1, \ldots, n-1, \\
\partial_t p^{z} &= - \partial_{z} H,\qquad\ \  p^{z,k}(1)= -\nabla \tilde U_k(z^{(k)}(1)) , \quad k=1, \ldots, n-1.
\end{split}
\end{equation*}
The computation of the differential system gives
\begin{multline*}
\prt_t p^{x,k}_i = - \sum_{j=1}^{m_k} \nabla_1(p_i^{x,k} \cdot K^{(k)}(x_i^{(k)}, x^{(k)}_j) \al^{(k)}_j) - \sum_{j=1}^{m_k} \nabla_1(\al_i^{(k)} \cdot K^{(k)}(x_i^{(k)}, x^{(k)}_j) p^{x,k}_j) \\
+ 2 \sum_{j=1}^{m_k} \nabla_1(\al_i^{(k)} \cdot K^{(k)}(x_i^{(k)}, x^{(k)}_j) \al^{(k)}_j) - (\la^{(k)} - \mu (x^{(k)} - z^{(k)})) ,
\end{multline*}
and
\begin{multline*}
\prt_t p^{z,k}_i = - \sum_{l=1}^{n-1}\sum_{j=1}^{m_l} \nabla_1(p_i^{z,k} \cdot K^{(n)}(z^{(k)}_i, z^{(l)}_j) \be^{(l)}_j) - \sum_{l=1}^{n-1}\sum_{j=1}^{m_l} \nabla_1(\be^{(k)}_i \cdot K^{(n)}(z^{(k)}_i, z^{(l)}_j) p^{z,l}_j) \\
+ 2 \sum_{l=1}^{n-1}\sum_{j=1}^{m_l} \nabla_1(\be^{(k)}_i \cdot K^{(n)}(z^{(k)}_i, z^{(l)}_j) \be^{(l)}_j) + \sum_{l=1}^n (\la^{(l)} - \mu (x^{(l)} - z^{(l)}))
\end{multline*}
The gradient of $L$ with respect to $(\al, \be)$ is then deduced from the partial differentials of $H$ with respect to these variables, yielding
\begin{equation*}
\begin{split}
\nabla_{\al^{(k)}} L &= K^{(k)}(x^{(k)}) (\al^{(k)} - p^{x,k}) , \\
\nabla_{\be} L &= K^{(n)}(z) (\be - p^{z}) .
\end{split}
\end{equation*}

Alternatively, on may choose to use the gradient relative to the dot product on $V_1 \times \ldots \times V_n$, which is simply given by 
\begin{equation*}
\begin{split}
\nabla_{\al^{(k)}} L &= \al^{(k)} - p^{x,k} , \\
\nabla_{\be} L &= \be - p^{z} .
\end{split}
\end{equation*}
The latter choice is simpler, and generally more efficient numerically.

\subsection{Sliding Interface }

Assume that the parameter sets $M_k$ are vertices of pure oriented, simplicial complexes $T_k$ of dimension $r_k < d$ (we will however only provide implementation details for codimension $d-r_{k}=1$). We let $F_k$ denote the set of facets of the $k$-th complex. We also assume that  $T_{1}, \ldots, T_{n-1}$ are disjoint and that $T_{n}$ is their  union, $T_{n} = \bigcup_{k=1}^{n-1} T_{k}$. We also let $F = \bigcup_{k=1}^{n-1} F_{k}$ (disjoint union).

The associated shape space is formed by functions $q_k: M_k \to \mR^d$ such that $q_{k}(f)$ is not degenerate (i.e., has maximal dimension) for all $f\in F_{k}$.  Each object is allowed to slide against the background. We will write $x^{(k)} = q^{(k)}(M_{k})$, $k=1, \ldots, n-1$, and $z^{(k)} = q^{(n)}(M_{k})$, $z = (z^{(1)}, \ldots, z^{(n-1)})$, in accordance with our previous notation.  If $f\in F_{k}$ is a facet in $T_k\subset T_{n}$, we discretize \eqref{eq:sliding} into
\begin{equation}
\label{eq:sliding.disc}
N^{(n)}(f)\cdot \left(\sum_{j\in f} (u_k(z^{(k)}_{j}) - u_n(z^{(k)}_{j})) \right)= 0,
\end{equation}
where  $N^{(n)}(f)$ is a $d\times (d-r_k)$ matrix spanning the normal space to $q^{(n)}(f)$, assumed to  be defined as a smooth function of $q^{(n)}$. If $r_k = d-1$, this is always possible, since $N^{(n)}$ is a vector that can be taken as the cross product of $z_{f,2} - z_{f,1}, \ldots, z_{f,d} - z_{f,1}$ where $z_{f,1}, \ldots, z_{f,d}$ is any labeling of the vertices of $q^{(n)}(f)$ ordered consistently with the orientation. 


\medskip

We now restrict to this case, with $d=3$, so that shapes are triangulated surfaces in $\mathbb R^{3}$, as discussed in Section \ref{sec:lddmm}. For $f\in F_{k}$ and $j\in f$, we denote by $e_{j,f}$ the edge $(z_{j''}^{(k)}-z_{j'}^{(k)})$ where $j'$ and $j''$ are the other two vertices of $f$ such that $(j,j',j'')$ is positively oriented.   Similarly, let $e'_{j,f} = (z_{j'}^{(k)} - z_{j}^{(k)})$ and $e''_{j,f} = (z_{j''}^{(k)} - z_{j}^{(k)})$ be the two edges stemming from $z_{j}^{(k)}$ so that 
\[
 e'_{j,f} \times e''_{j,f} = 2\,\mathrm{area}(q^{(n)}(f)) N^{(n)}(f) =: \tilde N^{(n)}(f)
\]
is the area-weighted positively oriented normal to $f$ in $q^{(n)}(M_{k})$. Note that $e_{j,f} = e''_{j,f}-e'_{j,f}$.

With this notation, we can rewrite the constraint in the form
\[
 \sum_{j\in f} \mathrm{det}(e'_{j,f}, e''_{j,f}, u_{k}(z^{(k)}_j) - u_{n}(z^{(k)}_j)) = 0.
\]
holding for all $f\in F_k$ and $k=1, \ldots, n-1$. 

Introducing a Lagrange multiplier $\lambda_{f}$ for each of these constraints, after reduction of the vector fields, which proceeds similarly to the identity constraints case, the augmented Lagrangian takes the form
\begin{multline*} L(\alpha,\beta, x, z) =
\frac12 \sum_{k=1}^{n-1} \int_0^1 \al^{(k)} \cdot (K^{(k)}(x^{(k)})\al^{(k)}) \, dt + \frac12 \int_0^1 \be\cdot (K^{(n)}(z)\be) \, dt +  \sum_{k=1}^{n-1} U_k(x^{(k)}(1))\\
+  \sum_{k=1}^{n-1} \tilde U_k(z^{(k)}(1)) - \sum_{k=1}^{n-1} \sum_{f\in F_{k}} \int_0^1 (\la_f \Ga^{(k)}_{f}(x^{(k)}, z) - \frac \mu 2 \Ga^{(k)}_{f}(x^{(k)}, z)^{2}) \, dt  ,
\end{multline*}
with
\[
\Gamma^{(k)}_{f}(x^{(k)}, z) := \sum_{j\in f} \mathrm{det}\left(e'_{j,f}, e''_{j,f}, K^{(k)}(z^{(k)}_{j}, x^{(k)}) \alpha^{(k)} - K^{(n)}(z_{j}^{(k)}, z) \beta\right).
\]

\medskip

We now compute the evolution equations for the co-states, as done with identity constraints.
 For $f\in F_{k}$ and $i\in M_k$, we have 
\begin{equation}
\label{eq:der.gamma.x}
\prt_{x_i^{(k)}} \Ga^{(k)}_{f} = 
\sum_{j\in f} \nabla_1(\al_i^{(k)} \cdot K^{(k)}(x_i^{(k)}, z^{(k)}_j) \tilde N^{(n)}(f)).
\end{equation}
Denoting
\[
\de^{(k)} (f) := \sum_{j\in f} (u_{k}(z^{(k)}_j) - u_{n}(z^{(k)}_j)) ,
\]
if $i\in f\in F_{k}$, then
\begin{multline}
\label{eq:der.gamma.z}
\prt_{z_i^{(k)}} \Ga^{(k)}_{f} = -  e_{i,f} \times \de^{(k)}(f) - \sum_{l=1}^{n-1}\sum_{j=1}^{m_l} \nabla_1(\tilde N^{(n)}(f) \cdot K^{(n)}(z_i^{(k)}, z^{(l)}_j) \be^{(l)}_j)\\
+  \sum_{j=1}^{m_k} \nabla_1( \tilde N^{(n)}(f) \cdot K^{(k)}(z_i^{(k)}, x^{(k)}_j) \al^{(k)}_j)
 -  \sum_{l=1}^{n-1}\sum_{j=1}^{m_l} \nabla_1(\be^{(k)}_i\cdot K^{(n)}(z^{(k)}_i, z^{(l)}_j) \tilde N^{(n)}(f)).
\end{multline}

Let  $p^{x,1}, \ldots, p^{x,n-1}$ and $p^{z}= (p^{z,1}, \ldots, p^{z,n-1})$ be the co-states. Let $\ga^{(k)}_f = \la_f - \mu \Ga^{(k)}_{f}$. For $i\in M_{k}$, let $F_{k}(i) = \{f\in M_{k}: i\in f\}$. Then
\begin{multline*}
\prt_t p^{x,k}_i = - \sum_{j=1}^{N_k} \nabla_1(p_i^{x,k} \cdot K^{(k)}(x_i^{(k)}, x^{(k)}_j) \al^{(k)}_j) - \sum_{j=1}^{N_k} \nabla_1(\al_i^{(k)} \cdot K^{(k)}(x_i^{(k)}, x^{(k)}_j) p^{x,k}_j) \\
+ 2 \sum_{j=1}^{N_k} \nabla_1(\al_i^{(k)} \cdot K^{(k)}(x_i^{(k)}, x^{(k)}_j) \al^{(k)}_j) - \sum_{f\in F_{k}} \ga^{(k)}_{f} \prt_{x_i^{(k)}} \Ga^{(k)}_{f} ,
\end{multline*}
and
\begin{multline*}
\prt_t p^{z,k}_i = - \sum_{j=1}^{N} \nabla_1(p_i^{z,i} \cdot K^{(n)}(z^{(k)}_i, z_j) \be_j) - \sum_{j=1}^{N} \nabla_1(\be^{(k)}_i \cdot K^{(n)}(z^{(k)}_i, z_j) p^{z}_j) \\
+ 2 \sum_{j=1}^{N} \nabla_1(\be^{(k)}_i \cdot K^{(n)}(z^{(k)}_i, z_j) \be_j) - \sum_{f\in F_{k}} \ga^{(k)}_{f} \prt_{z_i^{(k)}} \Ga^{(k)}_{f} ,
\end{multline*}
where $\prt_{x_i^{(k)}} \Ga^{(k)}_{f}$ and $\prt_{z_i^{(k)}} \Ga^{(k)}_{f}$ are given by \eqref{eq:der.gamma.x} and \eqref{eq:der.gamma.z}.

\medskip

For $f\in F_{k}$ ($k=1, \ldots, n-1$), we have
\begin{equation*}
\begin{split}
\prt_{\alpha^{(k)}} \Ga^{(k)}_{f} &= \sum_{j\in f} K^{(k)}(x^{(k)}, z_{j}^{(k)}) \tilde N^{(n)}_{j}(f) , \\
\prt_{\beta} \Ga^{(k)}_{f} &= -\sum_{j\in f} K^{(n)}(z, z_{j}^{(k)}) \tilde N^{(n)}_{j}(f) ,
\end{split}
\end{equation*}
Letting $\theta_j^{(k)} = \sum_{f\in F_k: j\in f} \gamma_f^{(k)} N_j^{(n)}(f)$, 
the gradient of $L$  in $\al$ and in $ \be$ is then given by
\begin{equation*}
\begin{split}
\nabla_{\al^{(k)}} L &= K^{(k)}(x^{(k)}) (\al^{(k)} - p^{x,k}) - K^{(k)}(x^{(k)}, z) \theta , \\
\nabla_{\be} L &= K^{(n)}(z) (\be - p^{z}) + K^{(n)}(z,z)  \theta
\end{split}
\end{equation*}
or, taking the Hilbert gradient,
\begin{equation*}
\begin{split}
\nabla_{\al^{(k)}} L &=  \al^{(k)} - p^{x,k} - K^{(k)}(x^{(k)})^{-1} K^{(k)}(x^{(k)}, z) \theta , \\
\nabla_{\be} L &= \be - p^{z} + K^{(n)}(z,z)  \theta.
\end{split}
\end{equation*}
In spite of it requiring the inversion of a linear system in the first equation, we found the latter version preferable to the $L^2$ gradient in our experiments.

\subsection{Remarks}

\subsubsection*{Existence of constrained solutions} 
It is important to note that, according to \cite[Theorem 1]{arguillere2014shape}, there always exists at least one solution of \eqref{eq:mult.oc}-\eqref{eq:mult.oc1} satisfying the constraints.

\subsubsection*{Convergence to surfaces} 
A question naturally arising is whether our discrete approximation using triangulations converges to the smooth setting as  triangles get smaller and smaller.  More precisely, assume that smooth initial surfaces $S^k_\init = q^{(k)}_\init(M_k)$ are triangulated, with increasingly fine triangulations $q^{(k,\ell)}_\init: M_{k,\ell} \to \mathbb R^3$, $\ell = 1,2,\ldots$,  where $M_{k,\ell}$ labels the vertices of a simplicial complex $T_{k,\ell}$ whose faces are $F_{k,\ell}$. We discuss whether minimizers $(u_{k,\ell}, k=1, \ldots, n)$ of the discrete problems have a subsequence  that converges to a minimizer $(u_k, k=1, \ldots, n)$ of the limit problem. 

Assume that the following  condition holds for the sequence of triangulations: 
\begin{itemize}
\item[(i)] We assume that for all $k$ and $\ell$, and for every $f\in F_{k,\ell}$, there exists an embedding  $\psi_{k, \ell}^f: t_{k, \ell}^f \to S^{(k)}_\init = q^{(k)}_\init(M_k)$ (where $t_{k, \ell}^f$ is the interior of the triangle $q_{k,\ell}(f)$)  such that $(\psi_{k, \ell}^f(t_{k, \ell}^f), f\in F_{k,\ell})$ partitions  $S^{(k)}_\init$ up to a negligible set and $\max_f \|\psi_{k,\ell}^f-\mathrm{id}_{t_{k,\ell}}\|_{1, \infty} \to 0$ when $\ell\to \infty$.
\end{itemize}
This conditions ensure that data attachment terms like those described in section  \ref{sec:examples} computed at diffeomorphic transformations $\varphi_{k,\ell} \circ q^{(k,\ell)}_\init$ converge to the same term computed at  $\varphi_{k} \circ q^{(k)}_\init$ as soon as $\varphi^{(k,\ell)}$ converges to $\varphi^{(k)}$ in $C^1(\mathbb R^3)$. Given this, we sketch  the argument leading to the consistency of the discrete approximations.
\medskip

For identity constraints, one can use \cite[Proposition 5]{arguillere2014shape}, which proves that, if the triangulations are nested (every vertex at step $\ell$ lies on the limit surface and is also a vertex at step ${\ell+1}$), then one can extract, from a corresponding sequence of identity-constrained optimal vector fields, a subsequence that converges towards an identity-constrained solution of \eqref{eq:mult.oc}-\eqref{eq:mult.oc1}.
\medskip

For sliding constraints, one cannot directly apply \cite[Proposition 5]{arguillere2014shape}, because the constraints are not nested, even when the triangulations are. To obtain a consistent approximation, we need to relax the discrete problems. More precisely, let $t\mapsto(u_1(t),\dots,u_n(t))\in V_1\times\dots\times V_n$ be a minimizer of the continuous problem with sliding constraints, and let $(\varphi_1,\dots,\varphi_n)$ denote the corresponding flow with $q^{(k)}(t)=\varphi_k(t)\circ q^{(k)}_{init}$ the corresponding deformation of $q^{(k)}_{init}$, and ${N}^{(k)}(t,x)=d\varphi_k(t)^TN^{(k)}_{init}(x)$ a normal to $S_k(t)=q^{(k)}(t)(M_k)$ at $q^{(k)}(t,x)$. In particular we have, for every $k=1,\dots,n$, every $x\in M_k$, and almost every time $t$,
$$
{N}^{(k)}(t,x)\cdot (u_k(t,q^{(k)}(t,x))-u_n(t,q^{(k)}(t,x))=0.
$$
Moreover, as $\sum_{k=1}^n\Vert u_k(t)\Vert^2$ is constant, both $u_k(t)$ and $du_k(t)$ are $\alpha$-Lipschitz for some positive constant $\alpha$ that does not depend on $t$ or $k$.

Now let $q^{(k,\ell)}(t)=\varphi_k(t)\circ q^{(k,\ell)}_{init}$ be the corresponding deformation of the discretization at step $\ell$. Recall that $N^{(k,\ell)}(t,f)$ denotes the unit normal to the triangle $q^{(k,\ell)}(t,f)$. We will prove:

\begin{itemize}
\item[(ii)] The discretized deformations at step $\ell$ satisfy the following relaxed sliding constraints
\[
\left\vert N^{(k,\ell)}(f)\cdot\left(\sum_{j\in f}(u_{k}(z^{(k, \ell)}_j)-u_{n}(z^{(k,\ell)}_j) \right)\right\vert \leq \varepsilon_\ell
\]
at almost every $t$ and for every face $f$ in $F_{k,\ell}$, for a suitably chosen sequence $\varepsilon_\ell>0$ going to $0$ as $\ell$ goes to $\infty$.
\end{itemize}

\noindent
Indeed, fix a face $f$ and an integer $\ell$. Define
$$
\bar{N}^{(k,\ell)}(t,f)=\frac{d\varphi_k(t)^TN^{(k,\ell)}(0,f)}{\vert d\varphi_k(t)^TN^{(k,\ell)}(0,f)\vert}
$$ 
for every time $t$. Note that assumption (i) implies that 
$$
\vert {N}^{(k,\ell)}(t,f)-\bar{N}^{(k,\ell)}(t,f)\vert\leq \gamma_\ell,
$$ 
for some sequence $\gamma_\ell$, independent of $f$ and $t$, and going to 0 as $\ell$ goes to infinity. Assumption (i) also implies that there exists a sequence $\eta_{\ell}$, independent of $f$ and going to $0$ as $\ell\rightarrow\infty$ such that for $y\in\psi^f_{k,\ell}(f)$, and $q^{(k)}_{init}(x)=y$, we have
$$
\vert {N}^{(k)}_{init}(x)-N^{(k,\ell)}(0,f)\vert+\sum_{j\in f} \vert q^{(k)}_{init}(x)-z_j^{(k,\ell)}(0)\vert\leq \eta_\ell.
$$
Some triangle inequalities, Gronwall's lemma, and the fact that $u_k$ and $du_k$ are $\alpha-$Lipschitz then imply
$$
\left\vert N^{(k,\ell)}(t,f)\cdot\left(\sum_{j\in f}(u_{k}(t,z^{(k, \ell)}_j)-u_{n}(t,z^{(k,\ell)}_j) \right)\right\vert \leq \varepsilon_\ell,
$$
with $\varepsilon_\ell=\alpha\eta_\ell e^\alpha+\gamma_\ell$ going to 0 as $\ell$ goes to infinity.

Consequently, if $(u_{1,\ell},\dots,u_{k,\ell})$ is a sequence of minimizers of the discretized problem at step $\ell$ with relaxed sliding constraints
\[
\left\vert N^{(k,\ell)}(t,f)\cdot\left(\sum_{j\in f}(u_{k,\ell}(t,z^{(k, \ell)}_j(t))-u_{n,\ell}(t,z^{(k,\ell)}_j(t)) \right)\right\vert \leq \varepsilon_\ell,
\]
we see that the infimum limit over $\ell$ of the respective discretized costs of $(u_{1,\ell},\dots,u_{n,\ell})$ is smaller than or equal to the cost of a minimizer of the continuous problem with sliding constraints. So to prove that a limit point of that sequence is a minimizer of the cost for the continuous problem with sliding constraints, all we need is to check that any such limit point does satisfy the constraints.

So let $(u_{1,\ell}, \ldots, u_{n,\ell})$ be a sequence of minimizers of the relaxed discrete problem at step $\ell$\footnote{It is easy to prove that such minimizers exist using the same method as that of \cite{arguillere2014shape}, and replacing equality constraints with inequality constraints.} that weakly converges to 
$(u_1, \ldots, u_n)$ in $V_1\times \cdots\times V_n$ (which is true for at least one subsequence of any minimizing sequence), then the associated flows  $\varphi_{k,\ell}$ and their first two space derivatives converge uniformly in time and space to the flows $\varphi_k$ and their first two derivatives.

From this it is easy to see that any sequence approximating $q^{(k)}_\init(x)$ as in assumption (i) is such that $q^{(k,\ell)}(t, j_\ell)\to q^{(k)}(t, x)$ and $N^{(k,\ell)}(t, j_\ell) \to N^{(k)}(t,x)$ at all times. Moreover, such a minimizing sequence must be, like all geodesics, such that  $\sum_{k=1}^n \|u_{k, \ell}(t)\|_{V_k}^2$ is constant in time, and smaller than the cost function associated to, say, $u_k=0$ for all $k$. This implies that the vector fields $u_{k,\ell}(t)$ are continuous uniformly in $k, \ell$ and $t$, which, combined with the continuity of the evaluation functionals in an RKHS implies that $u_{k,\ell}(q^{(k,\ell)}(t, j_\ell)) \to u_k(q^{(k)}(t,x))$ for all times. Consequently, one easily checks that each $(u_{1,\ell}, \ldots, u_{n,\ell})$ satisfy a relaxed version of the continuous version of the constraints, with a precision that goes to 0 as $\ell$ goes to infinity. This finishes proving that the constraints are satisfied with exactitude at the limit. 

\subsubsection*{Kernel derivatives}
Expressions similar to $\nabla_1(n \cdot K(x,y) \al) $ appear at multiple times in the previous computation (for some vectors $n$ and $\alpha$). 
For radial kernels ($K(x, y) = G(|x-y|^2)\Id_{\mR^d}$), we have
\[
\nabla_1(n \cdot K(x,y) \al) = 2 G'(|x-y|^2) (n\cdot\al) (x-y),
\]
which (slightly) simplifies the expressions. 
%
%

\subsubsection*{Sliding Interface -- Alternate Version}
As discussed in Section \ref{sec:multi.shape}, the sliding constraint can also be handled by introducing a new state variable $N$ that tracks a  vector (or frame) normal to the interface via \eqref{eq:normal.evol}. In the discrete case, one can discretize this equation by introducing states $N(f), f\in F$, indexed by the facets of $M$, and evolving according to
\[
\prt_t N(f) = -\frac 1{|f|}\sum_{i\in f}du^{(n)}(z_i)^T N(f) ,
\]
where $|f|$ is the number of vertices in $f$. The sliding constraints are now expressed in terms of the state variables in a more direct way, but with  a new co-state variable for the normals, bringing in an extra degree of complexity and increasing the computational cost. Note that the finite-dimensional reduction is still possible in this case, so that
$u^{(n)}(\cdot) = K^{(n)}(\cdot, z) \be$, and the evolution of the normals can be expressed in a form involving the differential of the kernel. This yields an adjoint system involving second derivatives of the kernel. We will not detail the computations in this paper, since they follow the same pattern as the other two that were already discussed (see \cite{sommer2013higher} for more examples on how higher-order variables can be handled in similar contexts). Note that this alternate version of the sliding constraints is slightly more general than the one discussed in the previous section, since it does not require a definition of a normal field that smoothly depends on the manifolds. 

\section{Experimental Results}\label{sec:exp}
\subsection{Synthetic Example}
\begin{figure}[htbp]
\centering
\includegraphics[trim=0 0.5in 0 0.5in, clip, width=0.49\textwidth]{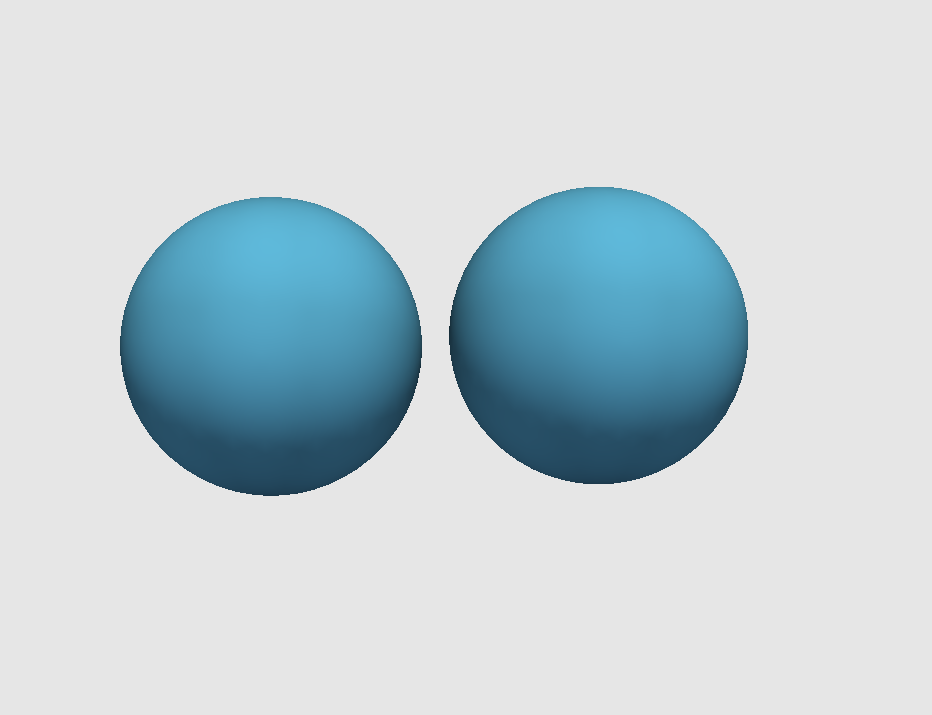}
\includegraphics[trim=0 0.5in 0 0.5in, clip, width=0.49\textwidth]{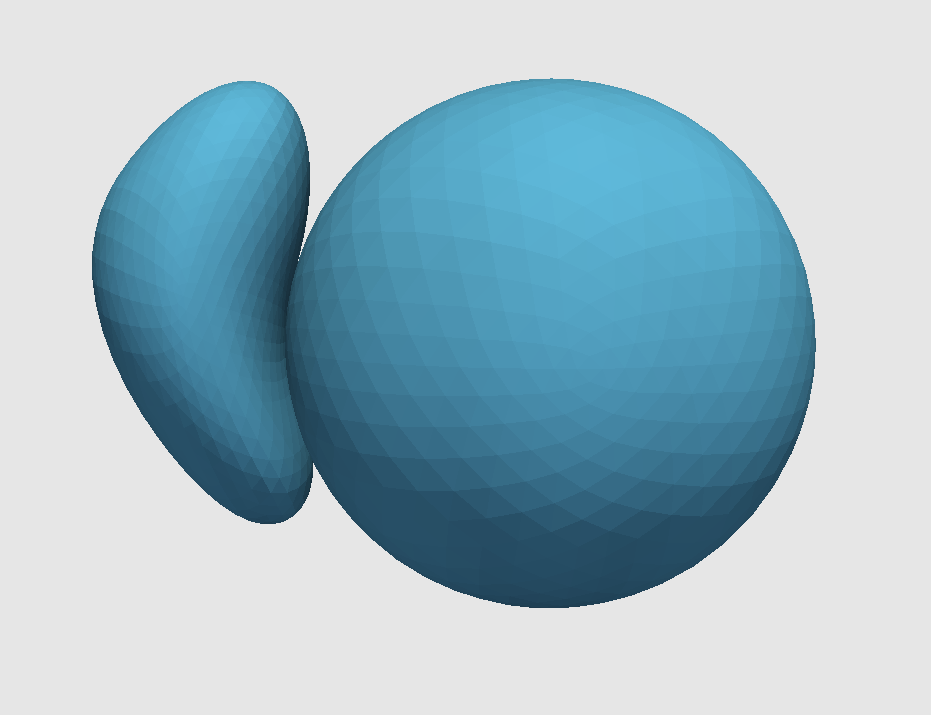}
\caption{\label{fig:twoBall} Template and target shapes for synthetic example}
\end{figure}

The first example is described in Figure \ref{fig:twoBall}. 
In this synthetic example, the template has two identical balls initially close to each other. In the target, the first ball (referred to as ``Ball A'') gets bigger, and ``impacts'' the other one (referred to as ``Ball B''), which assumes an oblong, non-convex shape (the target shapes slightly overlap, so that an exact homeomorphic match cannot be achieved). 

\begin{figure}[htbp]
\centering
\includegraphics[width=0.45\textwidth]{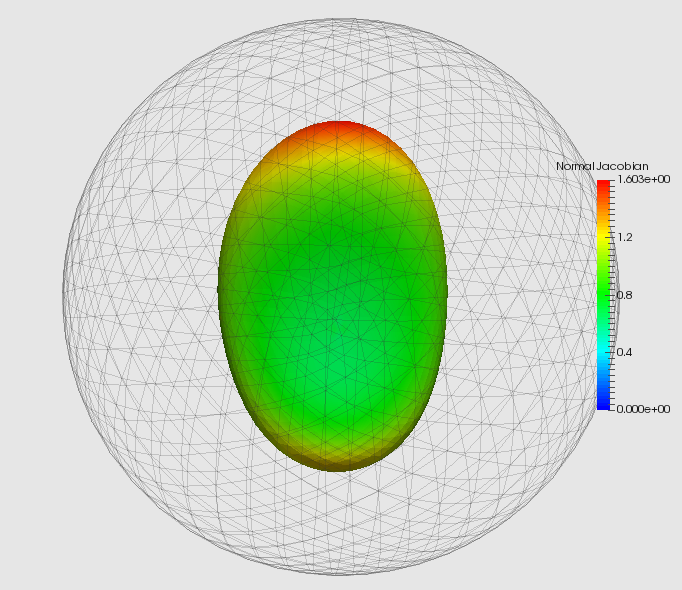}
\includegraphics[width=0.45\textwidth]{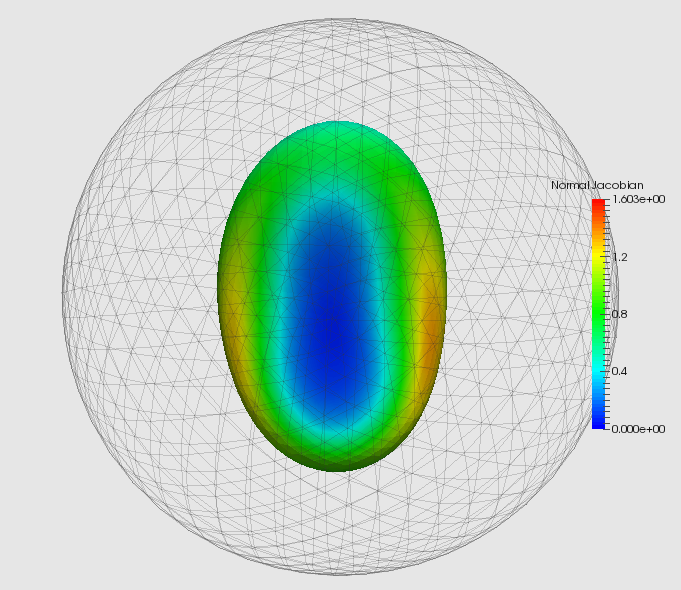}\\
\includegraphics[width=0.45\textwidth]{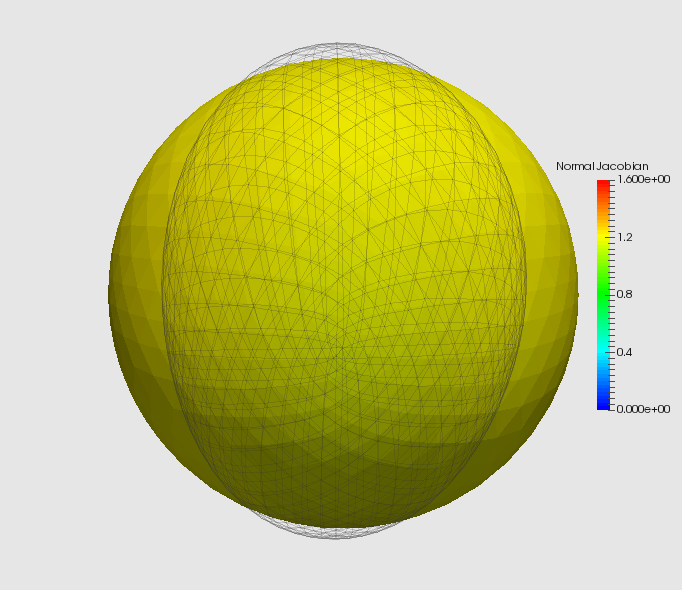}
\includegraphics[width=0.45\textwidth]{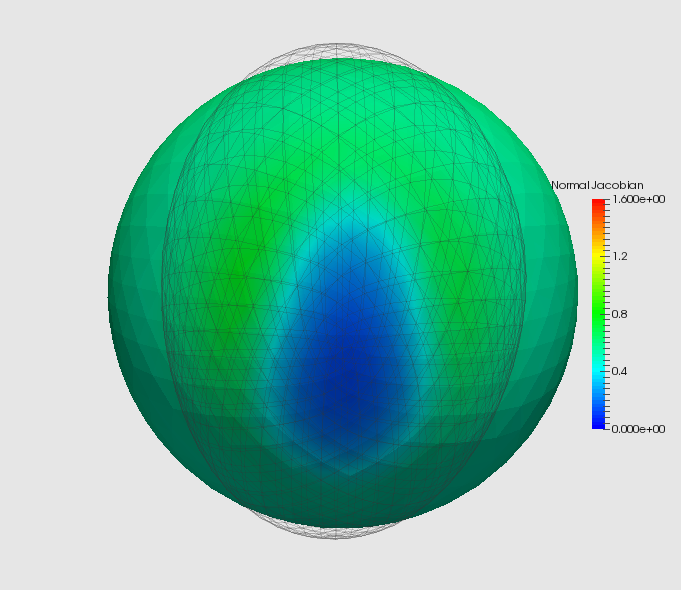}\\
\includegraphics[width=0.45\textwidth]{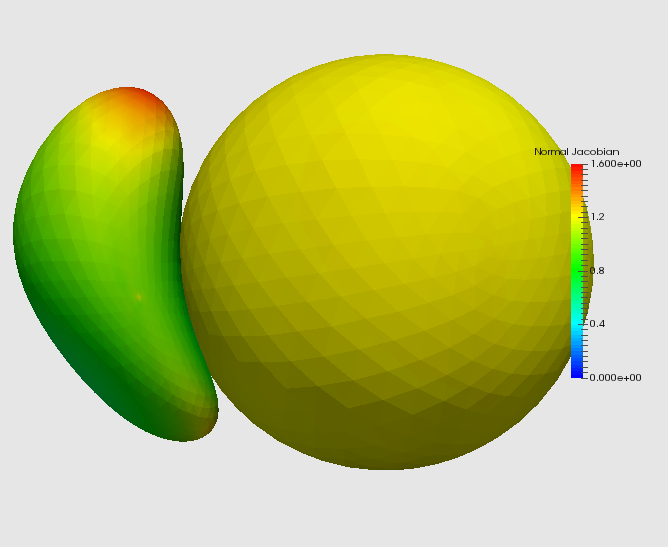}
\includegraphics[width=0.45\textwidth]{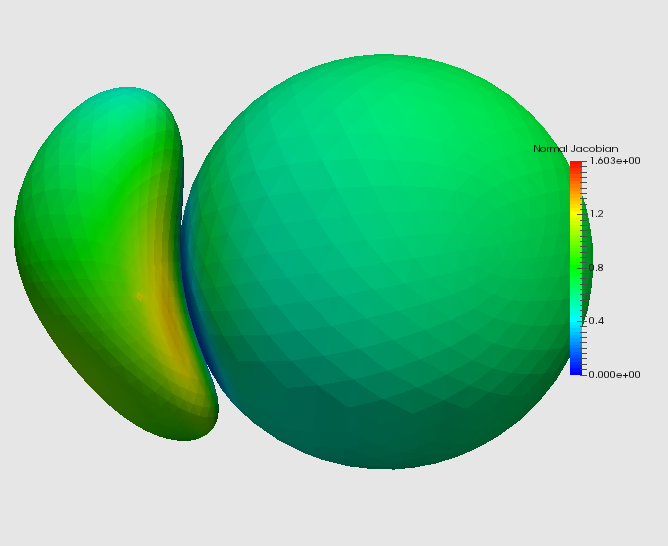}
\caption{\label{fig:twoBallStitchedNJ} Three views of the normal Jacobian with identity constraints: shape diffeomorphisms (left) and background diffeomorphism (right).}
\end{figure}

\begin{figure}[htbp]
\centering
\includegraphics[width=0.45\textwidth]{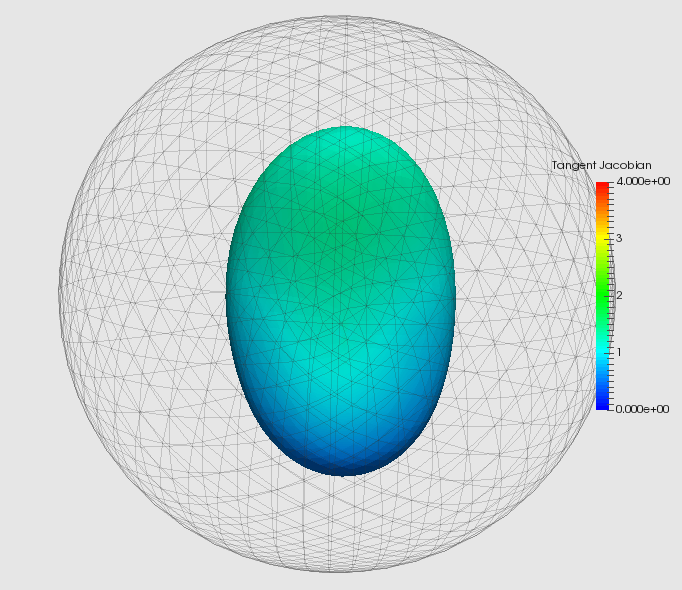}
\includegraphics[width=0.45\textwidth]{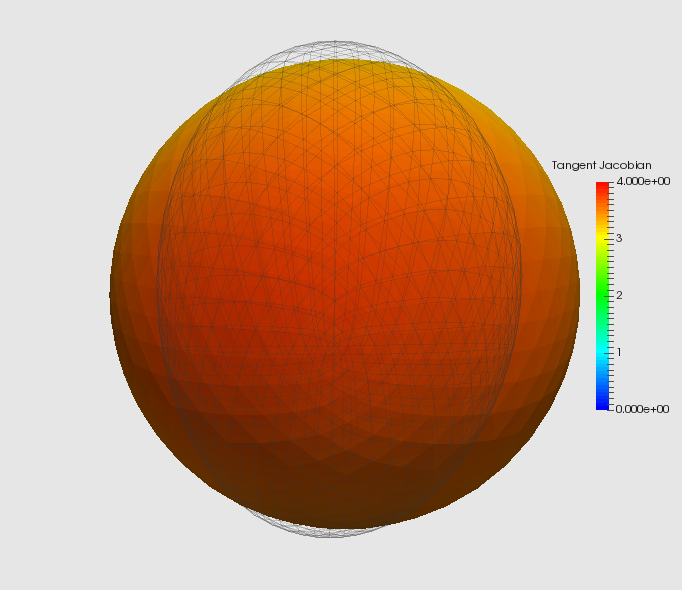}
\includegraphics[width=0.45\textwidth]{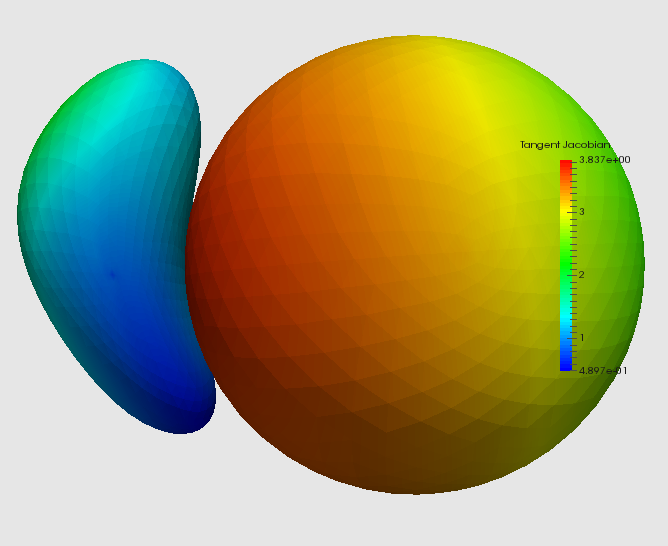}
\caption{\label{fig:twoBallStitchedTJin} Tangential Jacobian: shape and background diffeomorphisms (identity constraints).}
\end{figure}

\begin{figure}[htbp]
\centering
\includegraphics[width=0.45\textwidth]{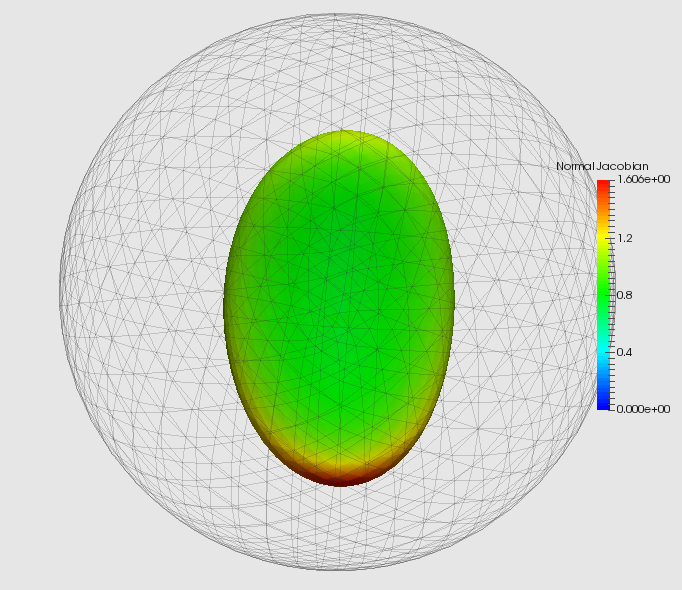}
\includegraphics[width=0.45\textwidth]{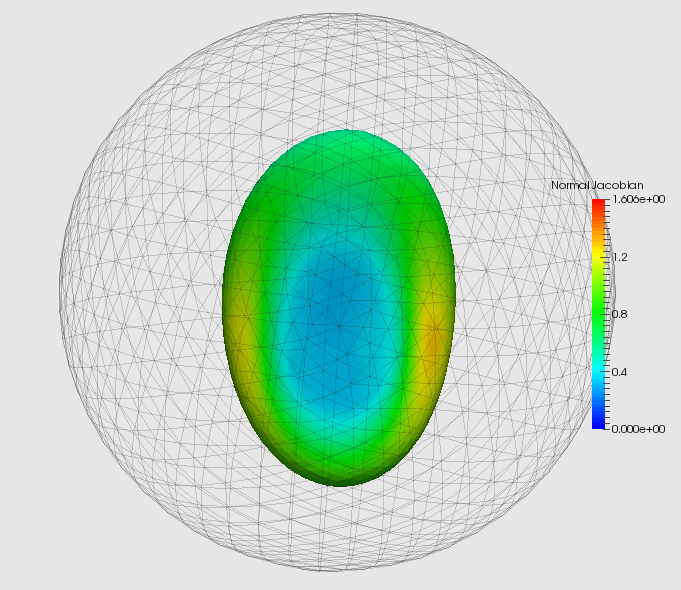}\\
\includegraphics[width=0.45\textwidth]{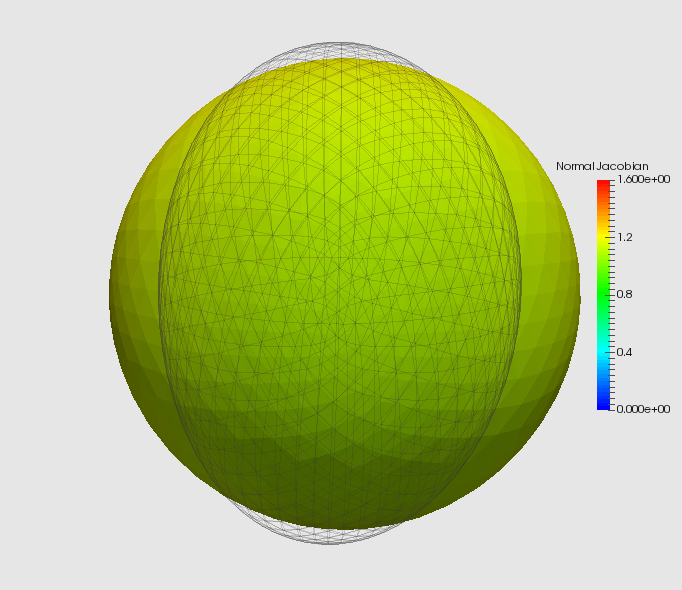}
\includegraphics[width=0.45\textwidth]{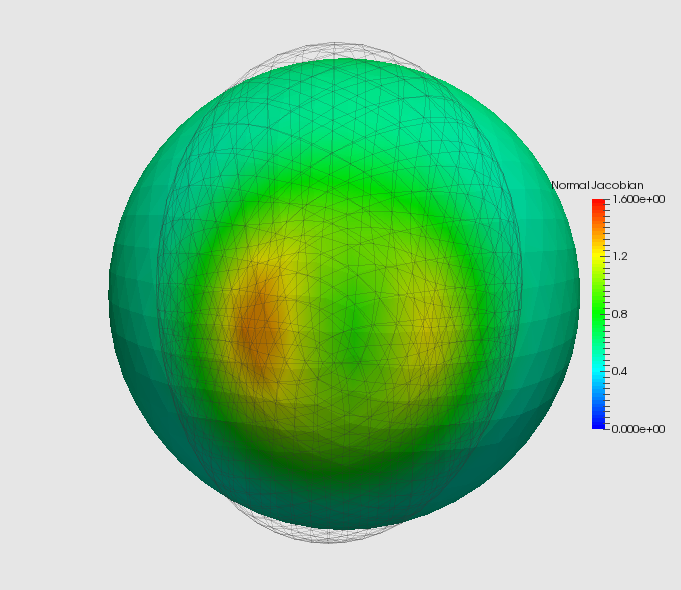}\\
\includegraphics[width=0.45\textwidth]{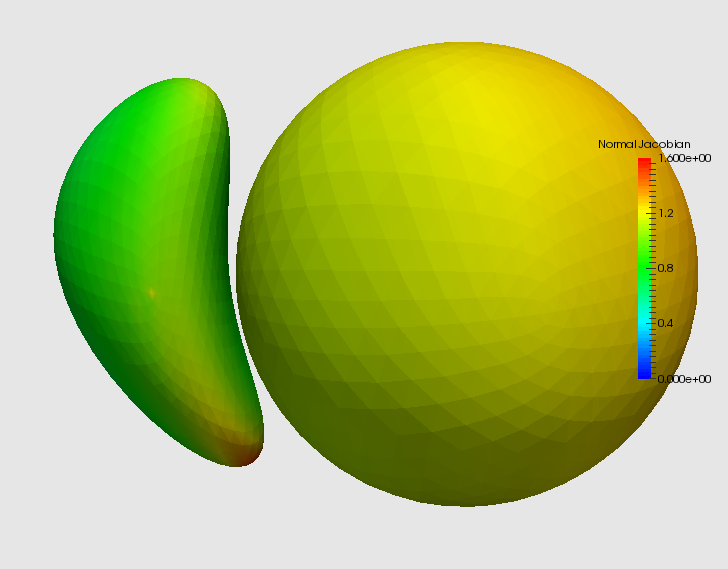}
\includegraphics[width=0.45\textwidth]{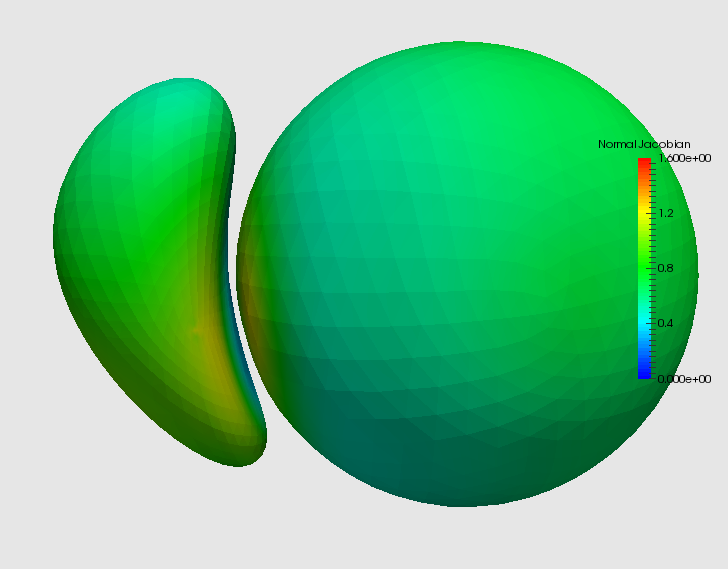}
\caption{\label{fig:twoBallSlidingNJin} Three views of the normal Jacobian with sliding constraints: shape diffeomorphisms (left) and background diffeomorphism (right).}
\end{figure}

\begin{figure}[htbp]
\centering
\includegraphics[width=0.45\textwidth]{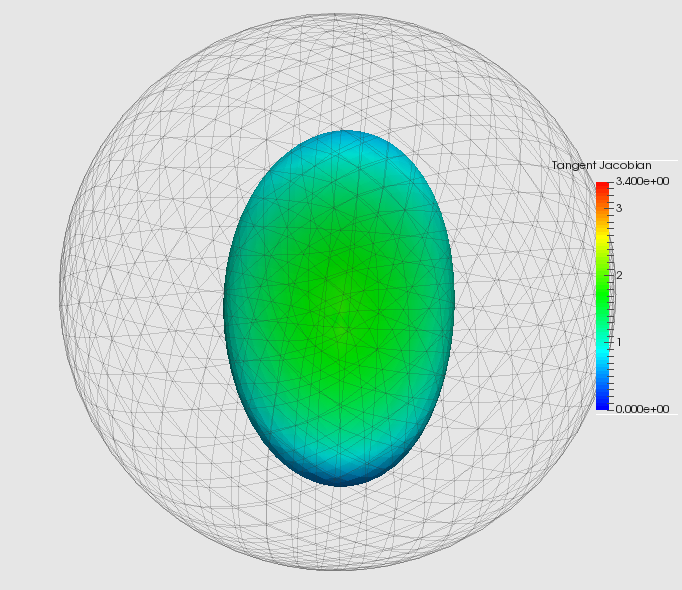}
\includegraphics[width=0.45\textwidth]{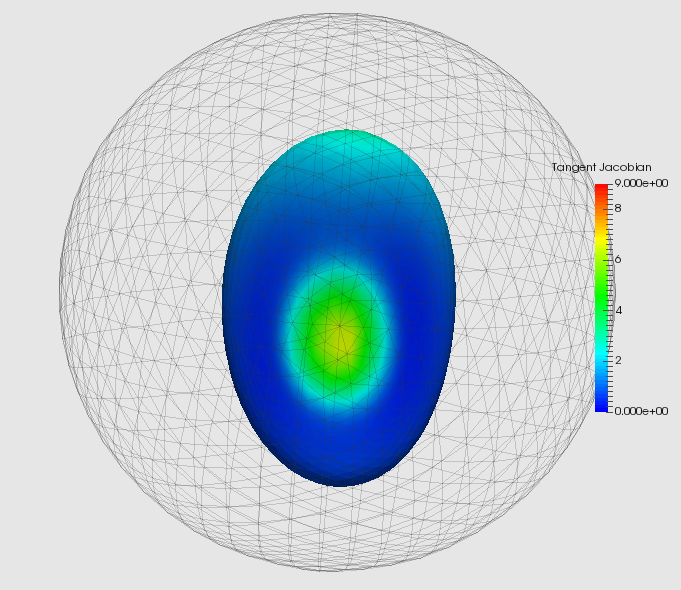}\\
\includegraphics[width=0.45\textwidth]{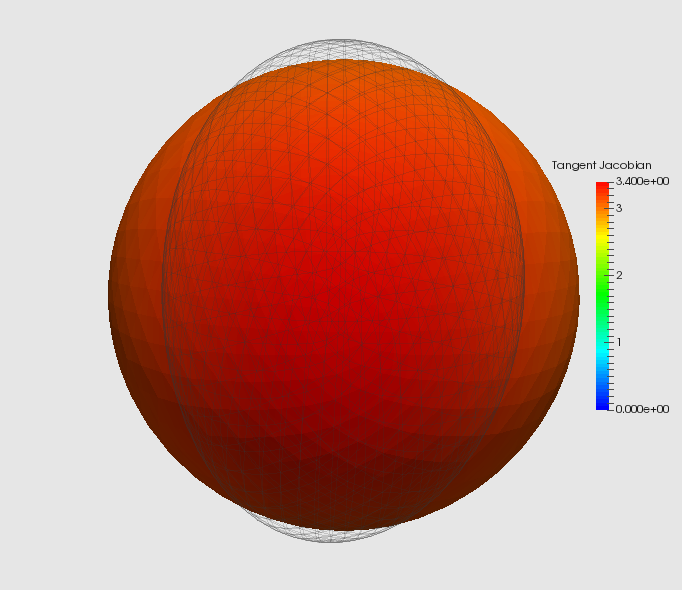}
\includegraphics[width=0.45\textwidth]{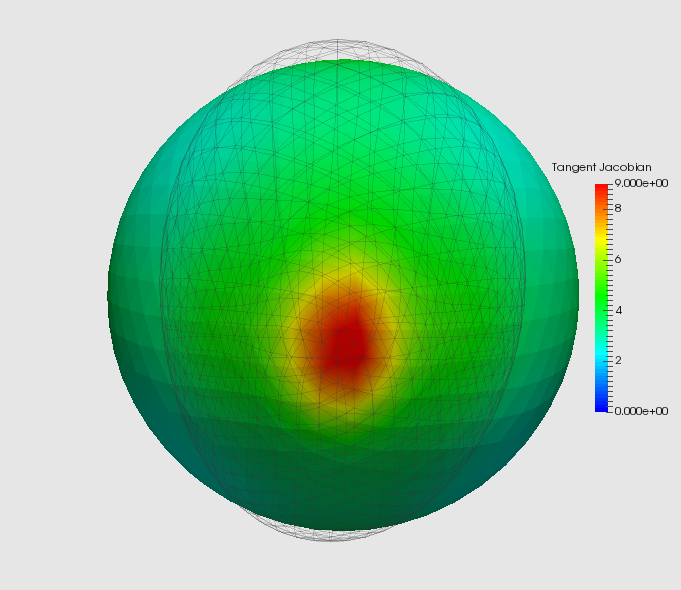}\\
\includegraphics[trim=4cm 0cm 0cm 0cm,clip,width=0.45\textwidth]{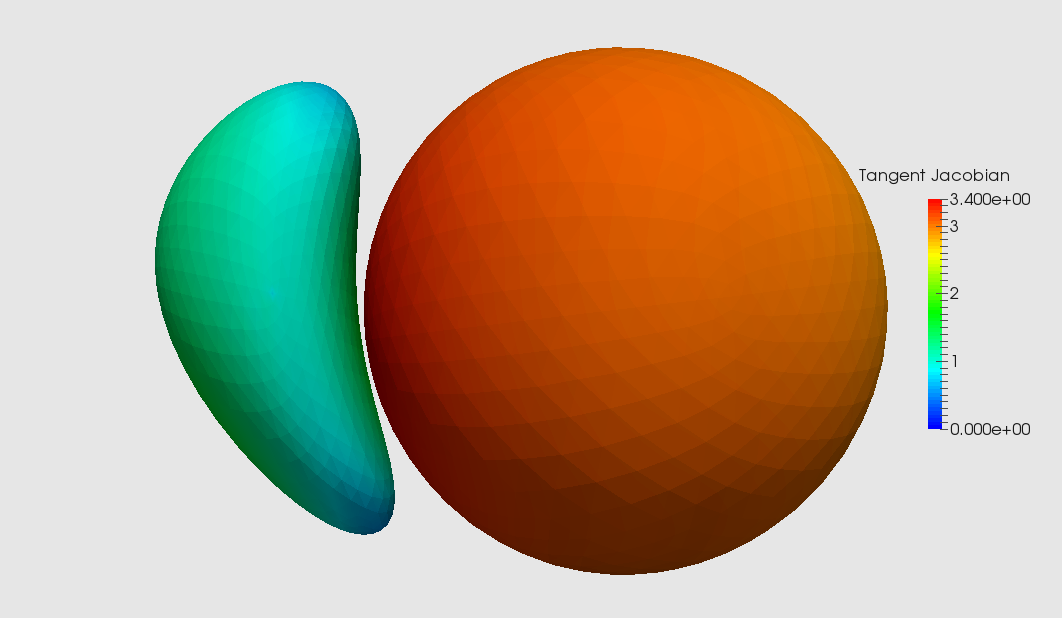}
\includegraphics[trim=4cm 0cm 0cm 0cm,clip,width=0.45\textwidth]{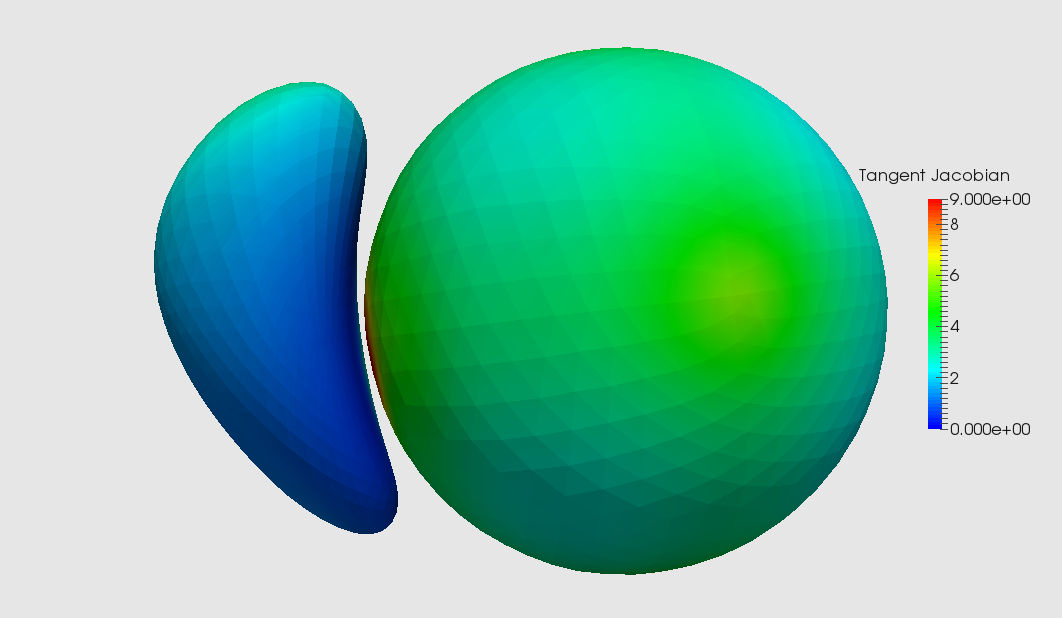}
\caption{\label{fig:twoBallSlidingTJin} Three views of the tangent Jacobian with sliding constraints: shape diffeomorphisms (left) and background diffeomorphism (right).}
\end{figure}

Our results, provided in Figures \ref{fig:twoBallStitchedNJ} to \ref{fig:twoBallSlidingTJin}, illustrate our multishape deformation method, and use two complementary deformation indexes:
\begin{itemize}
\item[(i)] the tangent Jacobian, which is the Jacobian determinant of the surface-to-surface transformations, and which measures the ratio between the areas of elementary surface patches at each point before and after deformation;
\item[(ii)] the normal Jacobian, which is the ratio of the Jacobian determinant (of the 3D diffeomorphism) to the tangent Jacobian, and which measures the ratio between the length of an infinitesimal line element normal to the surface after and before transformation. 
\end{itemize}
These indexes are mapped on the deformed template image, which is close to the target. 

Figure \ref{fig:twoBallStitchedNJ} compares the normal Jacobian of the shape and background deformations when using identity constraints. While shape diffeomorphisms characterize each shape transformation (uniform variation for Ball A, expansion at the top and compression otherwise for Ball B), the effect of compressing the space is clearly visible in the background deformation, when the two shapes get close to each other.  

Figure \ref{fig:twoBallStitchedTJin} provides the corresponding tangent Jacobian, which is identical for shape and background transformation since we are using identity constraints.

Figures \ref{fig:twoBallSlidingNJin} and \ref{fig:twoBallSlidingTJin} compare the normal and tangent Jacobians for the synthetic experiment with sliding constraints. Regarding the former, the most notable difference  is with Ball B, which shows an expansion pattern at the tips in its shape diffeomorphism  which is inverse of the one observed with identity constraints. One plausible explanation is that sliding constraints allow the two shapes to use translation-like motion to position themselves differently, without the need for limiting the amount of shear in the background that would have resulted from identity constraints. The second notable difference can be noted in the background diffeomorphism, in which compression is mostly observed with Ball B. In contrast with the identity constraints, the tangent Jacobians are very different between shape and background diffeomorphisms. Note that Figure \ref{fig:twoBallSlidingTJin} uses two different color scales for the left and right panels because of the strong difference between the ranges of the Jacobians in each case. The background deformation, in particular, has a huge tangent expansion around the impact location, which cannot be observed in the shape deformations. Note that both patterns in the sliding case are very different from the one that was observed in the identity case.

\begin{figure}[htbp]
\centering
\includegraphics[width=0.45\textwidth]{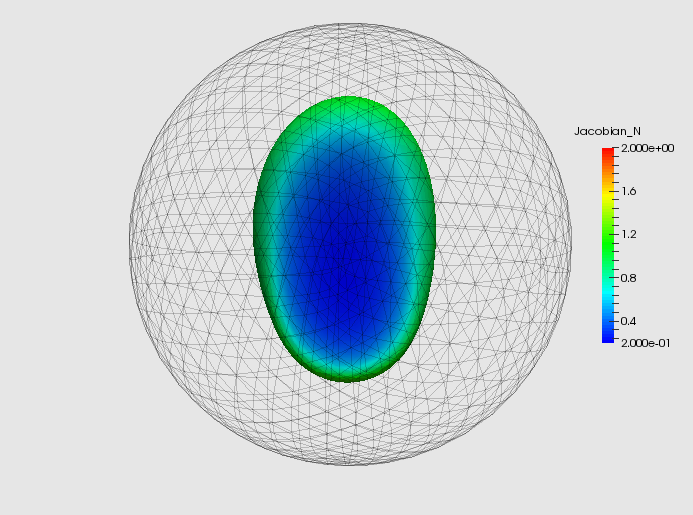}
\includegraphics[width=0.45\textwidth]{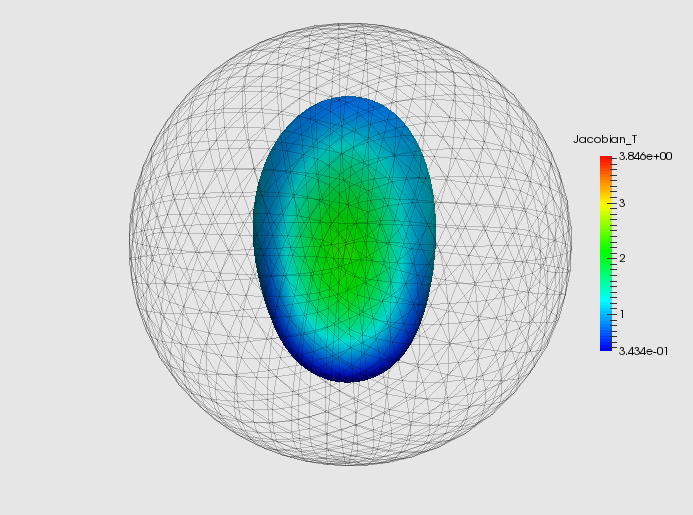}
\includegraphics[width=0.45\textwidth]{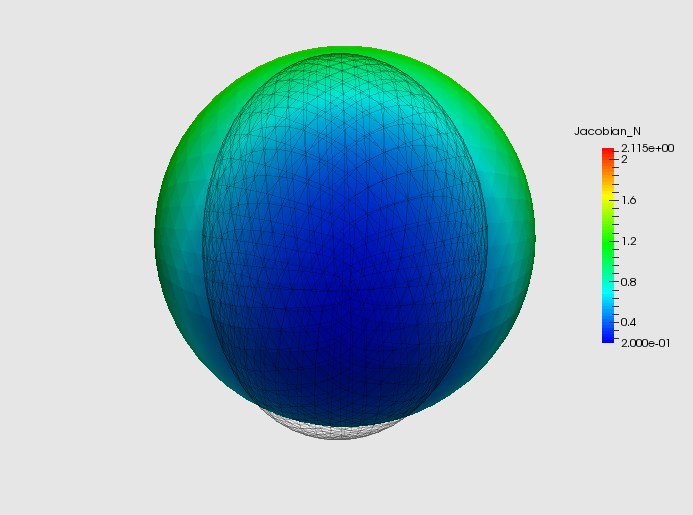}
\includegraphics[width=0.45\textwidth]{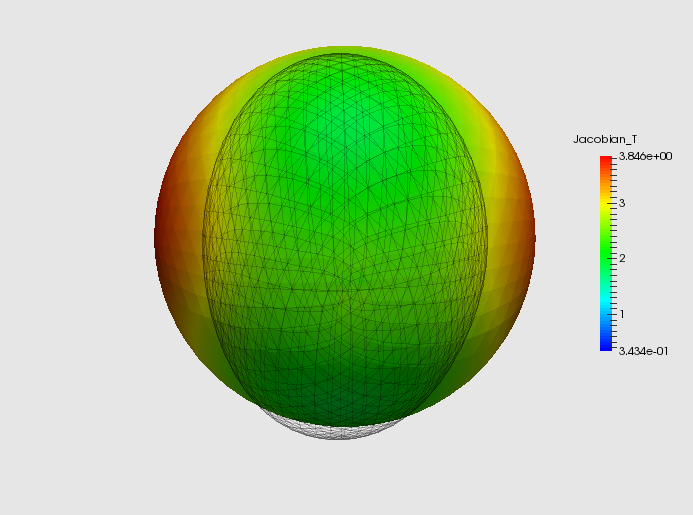}
\includegraphics[width=0.45\textwidth]{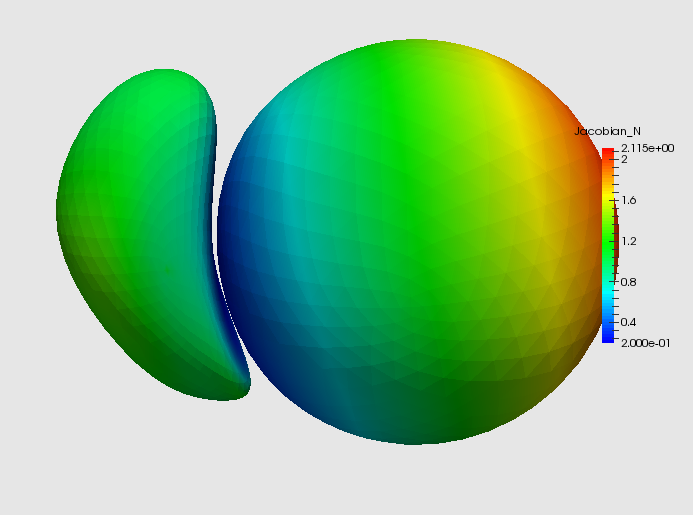}
\includegraphics[width=0.45\textwidth]{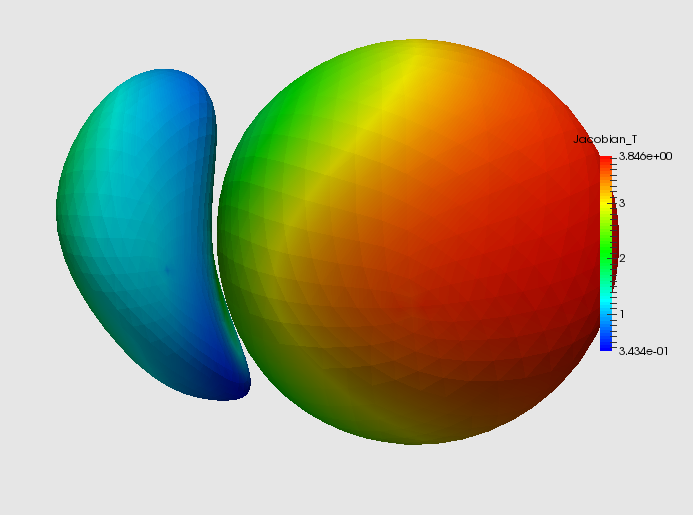}
\caption{\label{fig:twoBallsOneDiff} Three views of the Normal (left) and tangential Jacobians (right) when using a single diffeomorphism.}
\end{figure}

For comparison purposes, Figure \ref{fig:twoBallsOneDiff} provides the result of the LDDMM algorithm using a single diffeomorphism. One observes a very strong compression effect for the normal jacobian resulting in an expansion observed on the tangent jacobian on Ball B, that was not observed in any of the constrained examples. The nice uniform expansion in Ball A that could be observed in the sliding constraint case is not observed either.

\subsection{Subcortical Structures}
\begin{figure}[htbp]
\centering
\includegraphics[trim=0.75in 0.5in 0.75in 0.5in, clip, width=0.45\textwidth]{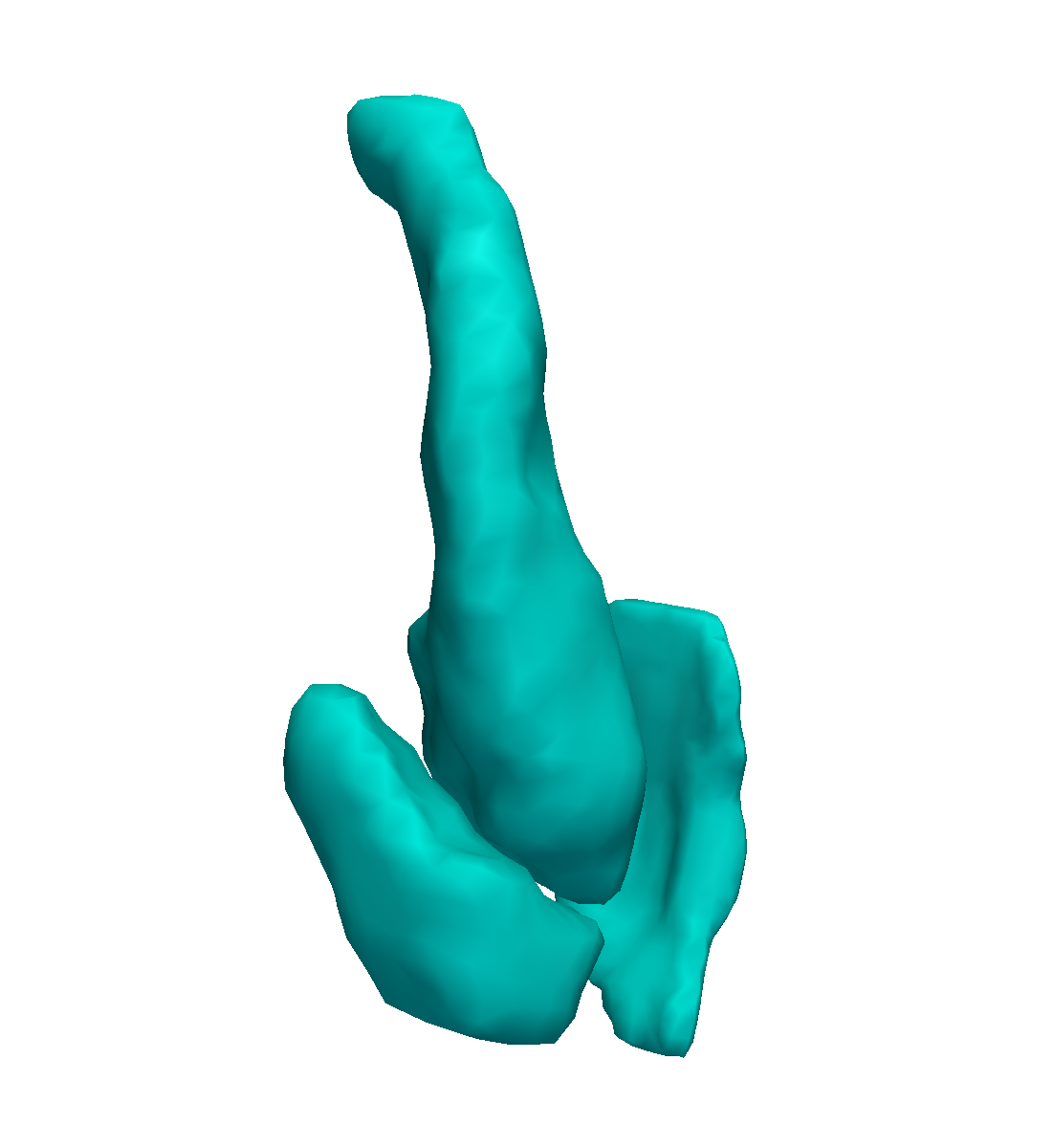}
\includegraphics[trim=0.75in 0.5in 0.75in 0.5in, clip, width=0.45\textwidth]{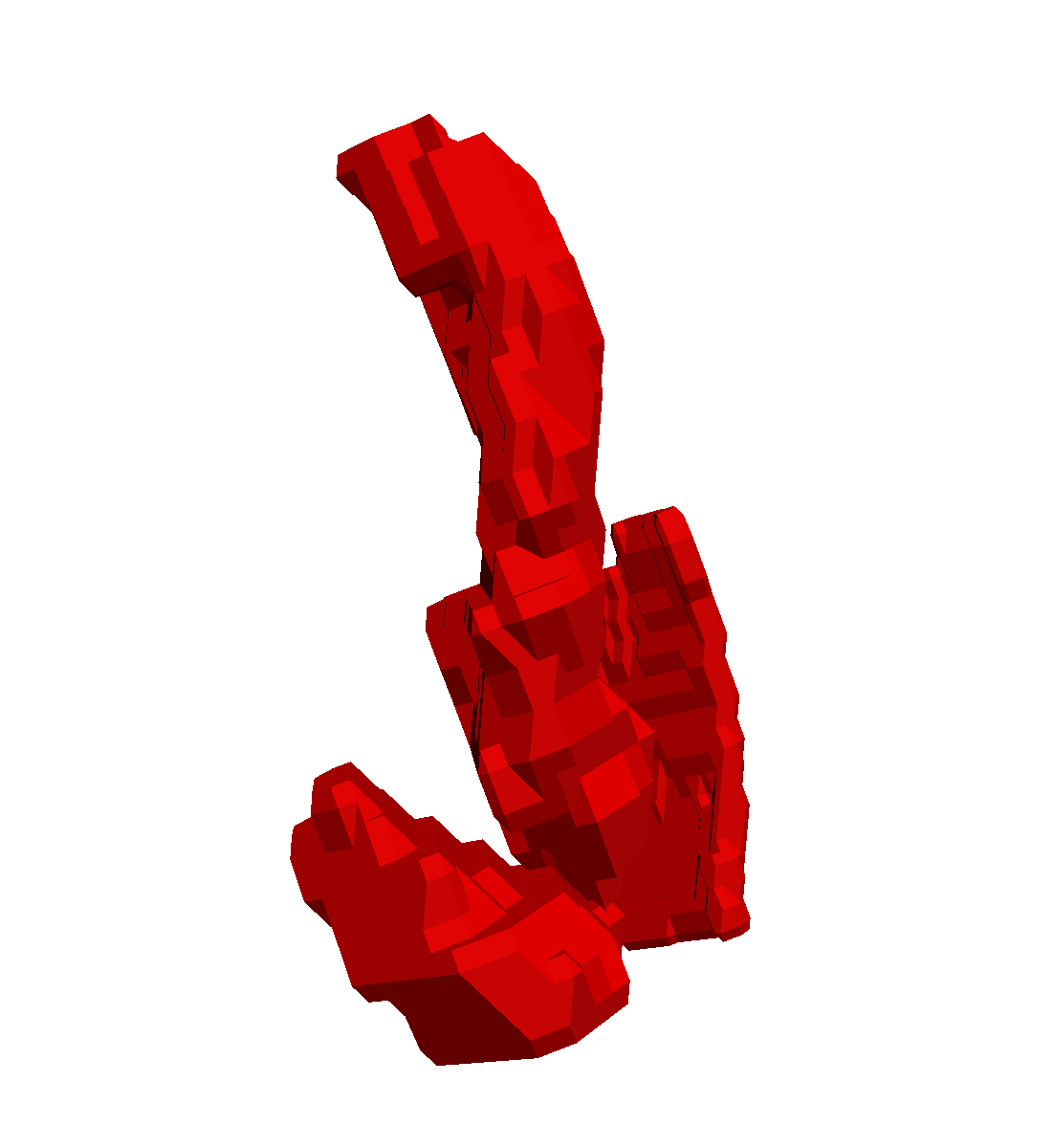}
\caption{\label{fig:biocard} Template (blue) and target (red) shapes for subcortical structures. The hippocampus is the central shape, with the amygdala on its left and the ERC on its right. }
\end{figure}

We now describe an example mapping a group of three subcortical structures: hippocampus, amygdala and entorhinal cortex (ERC). The template and target sets are represented in Figure \ref{fig:biocard}. 
One can observe shape changes in each structure, combined with a significant displacement of the ERC relative to the other two structures when comparing template to target. Because the structures were segmented independently, there is some overlap between the target hippocampus and amygdala. 

\begin{figure}[htbp]
\centering
\includegraphics[width=0.3\textwidth]{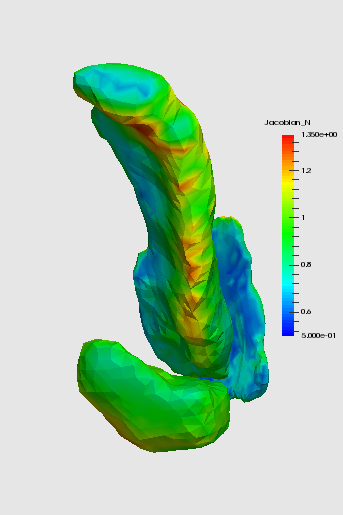}
\includegraphics[width=0.3\textwidth]{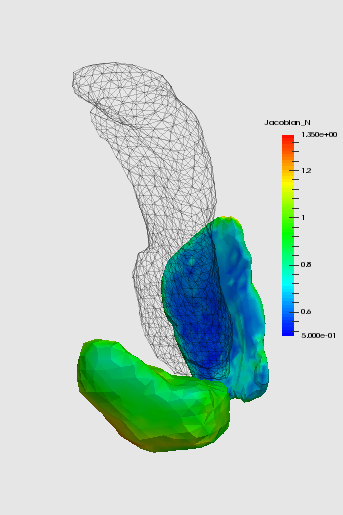}
\includegraphics[width=0.3\textwidth]{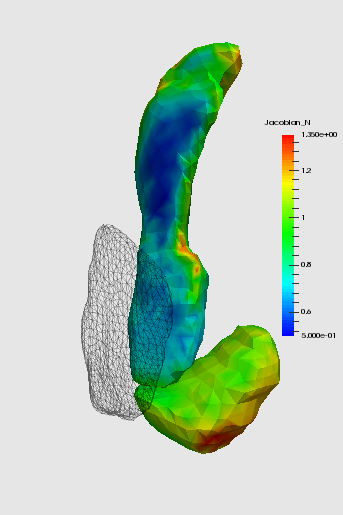}\\
\includegraphics[width=0.3\textwidth]{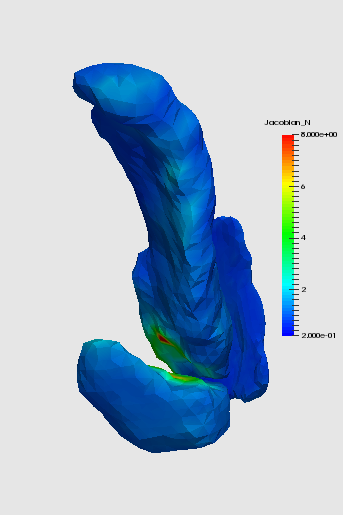}
\includegraphics[width=0.3\textwidth]{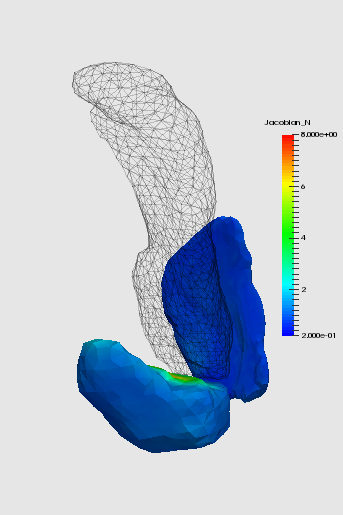}
\includegraphics[width=0.3\textwidth]{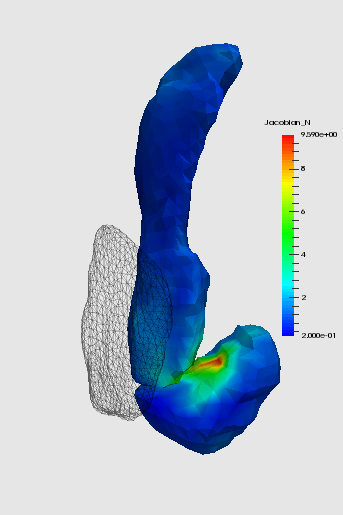}
\caption{\label{fig:biocardStitchedNJ} Three views of the normal Jacobian with identity constraints: First Row: Shape diffeomorphisms; Second Row: Background diffeomorphism.}
\end{figure}

\begin{figure}[htbp]
\centering
\includegraphics[width=0.3\textwidth]{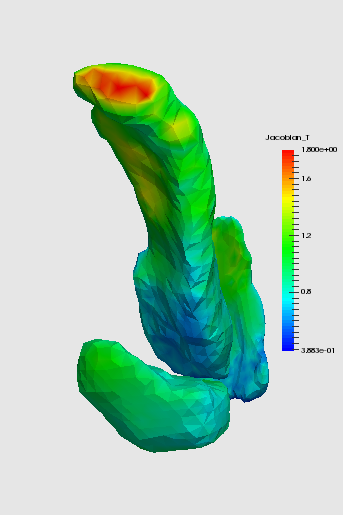}
\includegraphics[width=0.3\textwidth]{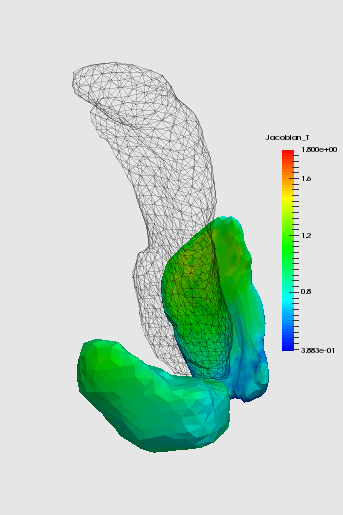}
\includegraphics[width=0.3\textwidth]{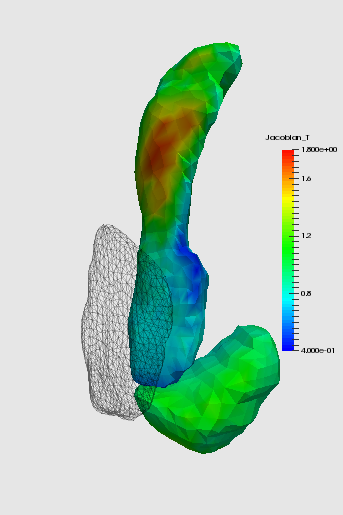}
\caption{\label{fig:biocardStitchedTJ} Three views of the tangent Jacobian with identity constraints: Shape and background diffeomorphisms.}
\end{figure}

\begin{figure}[htbp]
\centering
\includegraphics[width=0.3\textwidth]{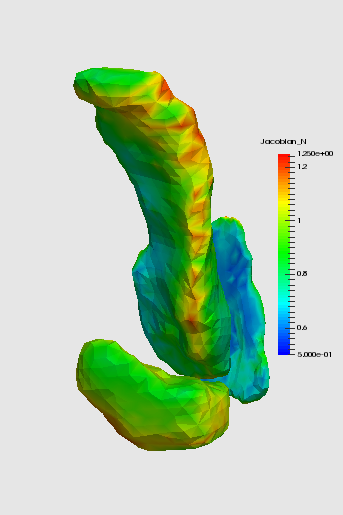}
\includegraphics[width=0.3\textwidth]{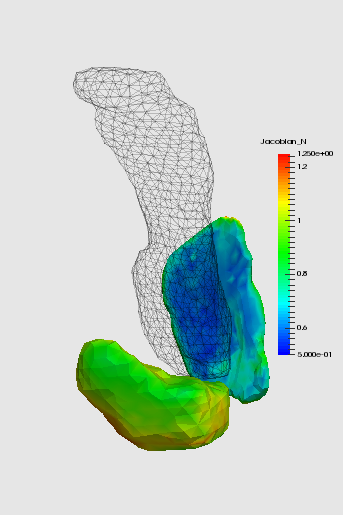}
\includegraphics[width=0.3\textwidth]{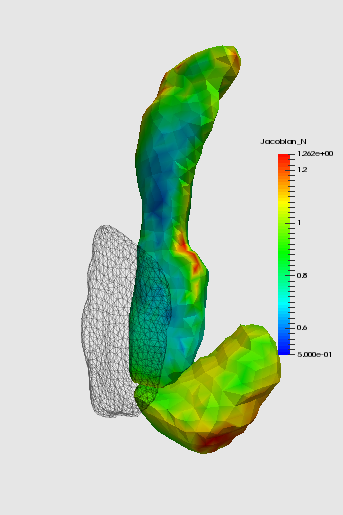}\\
\includegraphics[width=0.3\textwidth]{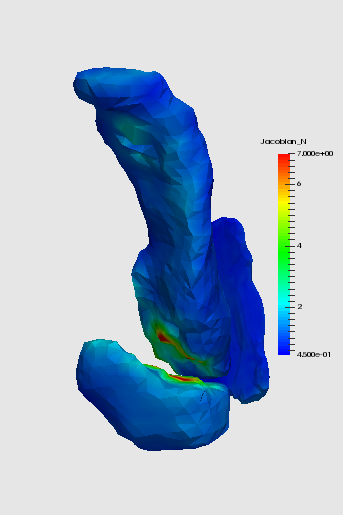}
\includegraphics[width=0.3\textwidth]{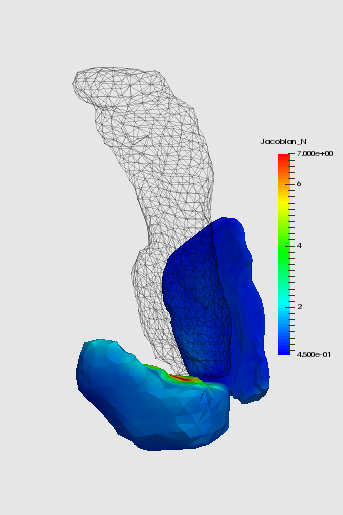}
\includegraphics[width=0.3\textwidth]{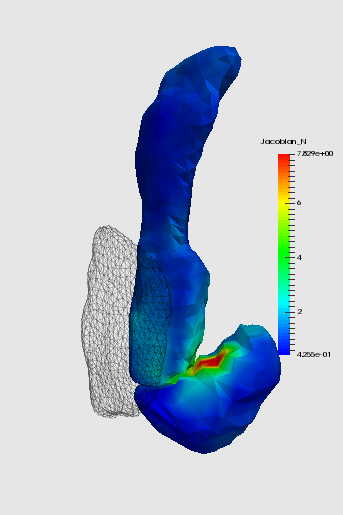}
\caption{\label{fig:biocardSlidingNJ} Three views of the normal Jacobian with sliding constraints: First Row: Shape diffeomorphisms; Second Row: Background diffeomorphism.}
\end{figure}

\begin{figure}[htbp]
\centering
\includegraphics[width=0.3\textwidth]{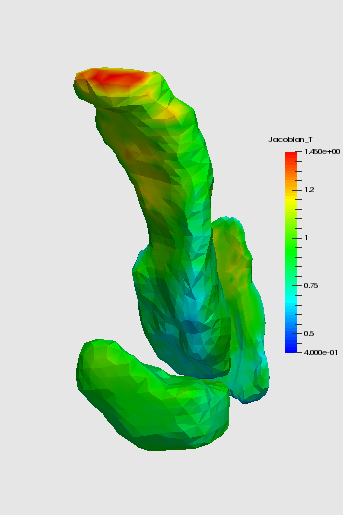}
\includegraphics[width=0.3\textwidth]{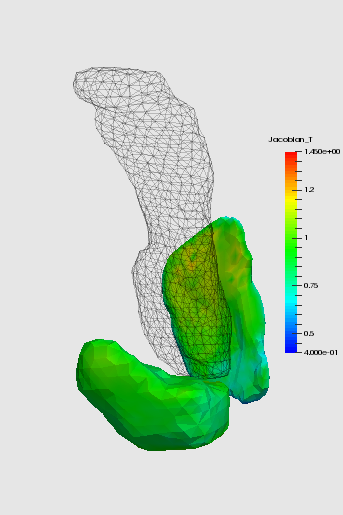}
\includegraphics[width=0.3\textwidth]{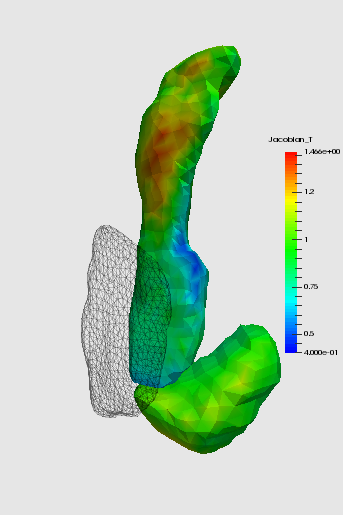}\\
\includegraphics[width=0.3\textwidth]{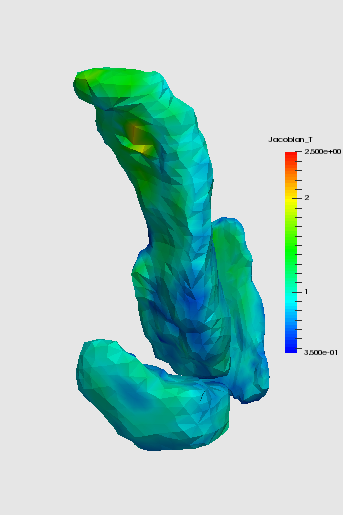}
\includegraphics[width=0.3\textwidth]{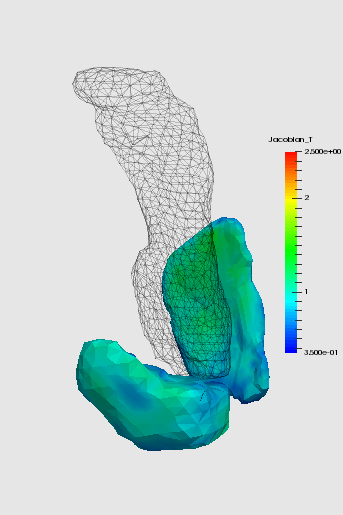}
\includegraphics[width=0.3\textwidth]{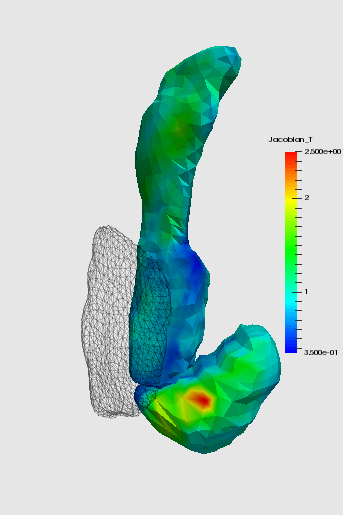}
\caption{\label{fig:biocardSlidingTJ} Three views of the tangent Jacobian with sliding constraints: First Row: Shape diffeomorphisms; Second Row: Background diffeomorphism.}
\end{figure}

\begin{figure}[htbp]
\centering
\includegraphics[width=0.3\textwidth]{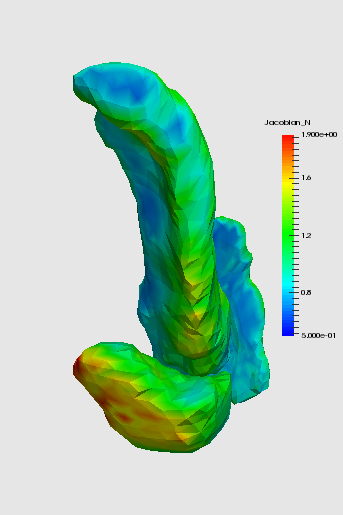}
\includegraphics[width=0.3\textwidth]{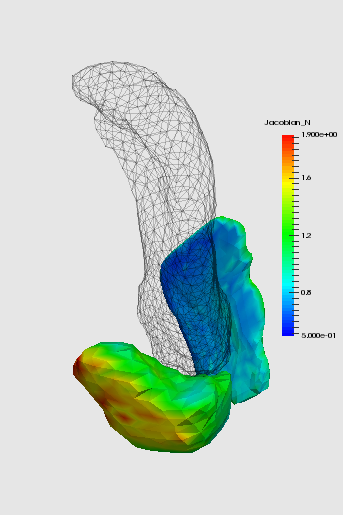}
\includegraphics[width=0.3\textwidth]{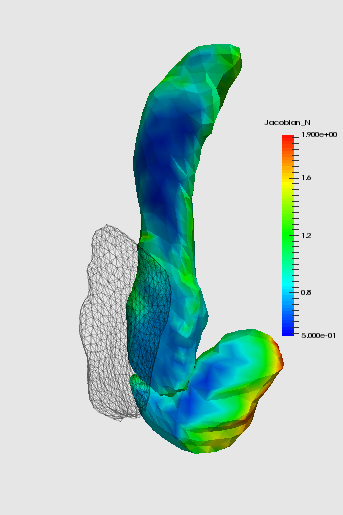}
\includegraphics[width=0.3\textwidth]{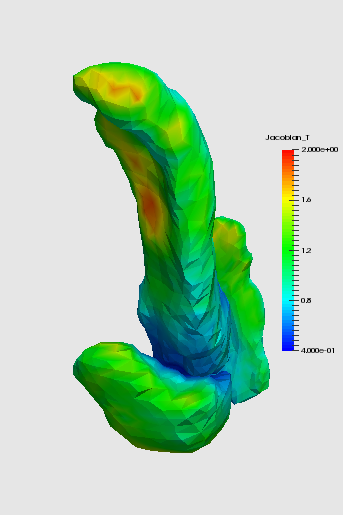}
\includegraphics[width=0.3\textwidth]{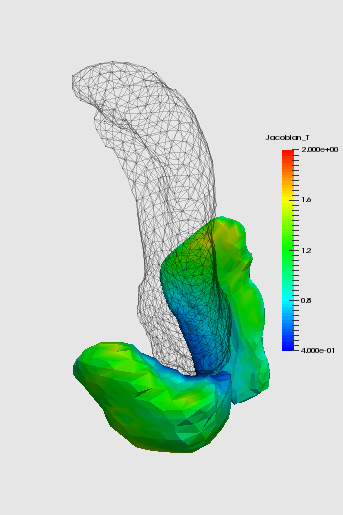}
\includegraphics[width=0.3\textwidth]{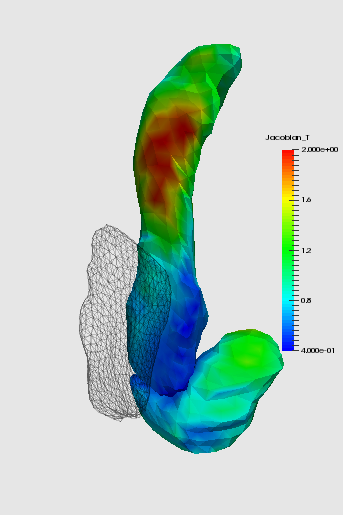}
\caption{\label{fig:biocardSingle} Three views of the normal (up) and tangent Jacobians (down) when using a single diffeomorphism.}
\end{figure}

Figures \ref{fig:biocardStitchedNJ} and \ref{fig:biocardStitchedTJ} provide the normal and tangent Jacobian obtained with identity constraints, while Figures \ref{fig:biocardSlidingNJ} and \ref{fig:biocardSlidingTJ} provide this information for sliding constraints. The two types of constraints provide similar deformation indices, especially for the normal jacobians (Figures \ref{fig:biocardStitchedNJ} and \ref{fig:biocardSlidingNJ}). Minor differences in the tangent jacobian can be observed  (Figures \ref{fig:biocardStitchedTJ} and \ref{fig:biocardSlidingTJ}). The deformation patterns associated to using a single diffeomorphism (Figure \ref{fig:biocardSingle}) are significantly different, though, exhibiting very strong compression, for example, where shapes are close to each other.

\begin{figure}[htbp]
\centering
\includegraphics[width=0.3\textwidth]{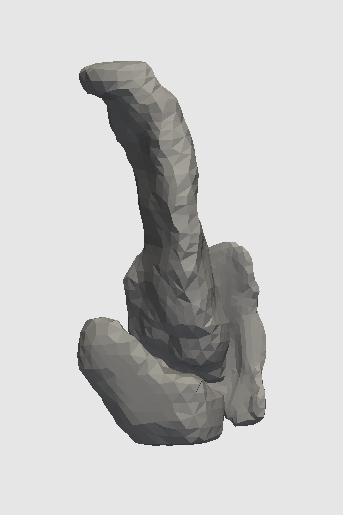}
\includegraphics[width=0.3\textwidth]{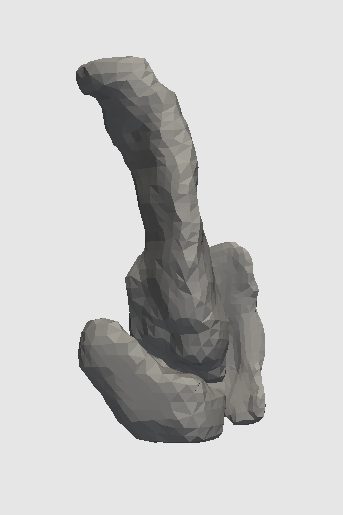}
\includegraphics[width=0.3\textwidth]{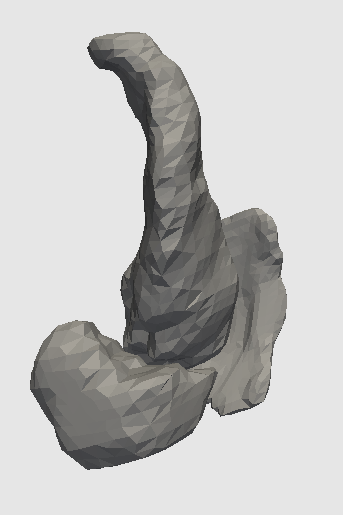}
\caption{\label{fig:biocardHalfway} Midpoint of the optimal deformation with multishape identity constraints (left), multishape sliding constraints (center) and single diffeomorphism (right).}
\end{figure}

%

\section{Discussion}
The previous approach provides a solution, using constrained optimal control, of the important issue of dealing with multiple objects with varying deformation properties for registration. We have focused on surface matching, numerically dealing with constraints using an augmented Lagrangian method. Note that a similar approach was introduced for plane curves in \cite{arguillere2014shape}. 

The formulation is quite general and can accommodate constraints in various forms, including the examples discussed in Section \ref{sec:multi.shape}. The investigation of these additional applications will be the subject of future work. One of the limitations of the present implementation is the slow convergence of the augmented Lagrangian procedure, for which each minimization step is, in addition, high dimensional and computationally demanding. One possible alternative can be based on solving the optimality conditions \eqref{eq:pmp.gen} (which hold in the discrete case) by means of a numerical shooting method. This approach has, however, its own numerical challenges, because solving \eqref{eq:pmp.gen} requires the determination of $\lambda$ such that the last equation (constraint) is satisfied, and this leads to a possibly ill-posed problem for systems in large dimension (see \cite{arguillere2014shape} for additional details).

We have illustrated our examples using deformation markers derived from the jacobian determinant. This markers are routinely used in shape analysis studies and led to important conclusion in computational anatomy. When dealing with multiple shapes, however, figures \ref{fig:twoBallsOneDiff} and \ref{fig:biocardSingle} show that, when using the classical LDDMM method with multiple shapes, these markers becomes as much, if not more, influenced by interactions between the shapes as by the changes in the shapes themselves. For this reason, multi-shape computational anatomy studies have applied registration methods separately to each shape, without ensuring that the obtained diffeomorphisms are consistent with each other. This limitation is addressed in the present paper, in which we exhibit deformation markers that are meaningful in describing tangential and normal surface stretching, while being consistently associated to a global transformation of the space.

\bibliographystyle{plain}
\bibliography{multishapes}

\end{document}